\documentclass{article}
\pdfoutput=1
\usepackage{amsmath}
\usepackage{amsthm}
\usepackage{amsfonts}
\usepackage{amssymb}
\usepackage{bm}
\usepackage{graphicx}
\usepackage{xspace}
\usepackage{xcolor}
\usepackage[colorlinks=true,citecolor=darkgreen]{hyperref}

\definecolor{darkgreen}{rgb}{0,0.6,0}

\hypersetup{pdftitle={The Bring sextic of equilateral pentagons}}
\hypersetup{pdfauthor={Lyle Ramshaw}}

\newenvironment{describe}{\begin{description}}{\end{description}}

\renewcommand\_{\hbox{-}}

\newcommand{\reals}{\mathbb{R}}
\newcommand{\ints}{\mathbb{Z}}
\newcommand{\complexes}{\mathbb{C}}

\DeclareMathOperator{\arccot}{arccot}

\DeclareMathOperator{\arccosh}{arccosh}

\DeclareMathOperator{\dev}{dev}
\DeclareMathOperator{\Span}{Span}

\mathchardef\ldotp="613A   % Take dots from math italic, not roman
\newcommand{\dd}{\mathbin{\ldotp\!\ldotp}}

\newcommand{\mystrut}[1]{\rule[0pt]{0pt}{#1}}
\newcommand{\mydepthstrut}[1]{\rule[-#1]{0pt}{#1}}

\newcommand{\edge}{\varepsilon}

\newcommand{\s}{\kern -0.3pt 's\xspace}
\renewcommand{\t}{\kern -0.8pt 't\xspace}

\newcommand{\Ltract}{$L$\_tract}
\newcommand{\Rtract}{$R$\_tract}

\newcommand{\LSR}{\texttt{LSR}}
\newcommand{\RSL}{\texttt{RSL}}

\newcommand{\cWise}{\texttt{CW}}
\newcommand{\ccWise}{\texttt{CCW}}

\newcommand{\smangle}{{\scriptscriptstyle\angle}}
\newcommand{\Bsextic}{\mathcal{B}}
\newcommand{\RiemannEfive}{\mathcal{E}_5}

\newcommand{\imagunit}{\mkern 1mu\mathbf{i}\mkern 1mu}
\newcommand{\subj}{_{\mkern -1mu j}}

\newcommand{\flect}{\text{\small$\updownarrow$}}

\newcommand{\cR}{\textsf{\small R}}
\newcommand{\cO}{\textsf{\small O}}
\newcommand{\cG}{\textsf{\small G}}
\newcommand{\cB}{\textsf{\small B}}
\newcommand{\cP}{\textsf{\small P}}

\begin{document}

\title{The Bring sextic 
of equilateral pentagons\footnote{Version of \today.}}
\author{Lyle Ramshaw\footnote{Email: \texttt{lyle.ramshaw@gmail.com}%
\quad\protect\includegraphics[scale=0.06]{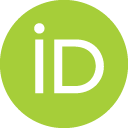}
\protect\href{https://orcid.org/0000-0002-9113-1473}{//orcid.org/0000-0002-9113-1473}} \ (retired)}
\date{\vspace{-12pt}}
\maketitle

\begin{abstract}
Consider equilateral pentagons $V_1\cdots V_5$ in the Euclidean plane.  When we identify pentagons that differ by translation, rotation, and magnification, the moduli space of possible shapes that we get is an oft-studied polygon space: a $2$\_manifold~$E_5$ known topologically to be a quadruple torus (genus~$4$).  We study~$E_5$ geometrically, our goal being a conformal map of that terrain of possible shapes.  The differential geometry that we use is all due to Gauss, though much of it is named after his student Riemann.
%, who generalized it to abstract $d$\_manifolds with $d\ge 2$.

The manifold $E_5$ inherits a Riemannian metric from the Grassmannian approach of Hausmann and Knutson, a metric $e_5$ under which $E_5$ has $240$ isometries: an optional reflection combined with any permutation of the order in which the five edge vectors $V_{k+1}-V_k$ get assembled into a pentagon.  Giving $E_5$ the conformal structure imposed by~$e_5$ yields a compact Riemann surface of genus~$4$ with $120$ automorphisms: the $120$ isometries that preserve orientation.  But there is only one Riemann surface with those properties: the \emph{Bring sextic}.   So $(E_5,e_5)$ conformally embeds into the hyperbolic plane, like the Bring sextic, as a repeating pattern of $240$ triangles, each with vertex angles of $\frac{\pi}{2}$, $\frac{\pi}{4}$, and $\frac{\pi}{5}$.  That conformal map realizes our goal.

To plot pentagons on our map, we compute an initial pair of isothermal coordinates for~$E_5$ by solving the Beltrami equation \`a la Gauss.  We then use a conformal mapping to convert one of those isothermal triangular regions into a Poincar\'e projection of a $(\frac{\pi}{2},\frac{\pi}{4},\frac{\pi}{5})$ hyperbolic triangle.
\end{abstract}


\vspace{5 pt}

\section{Introduction}
\label{sect:Intro}

\begin{figure}
	\begin{center}
	    \includegraphics[scale=0.6]{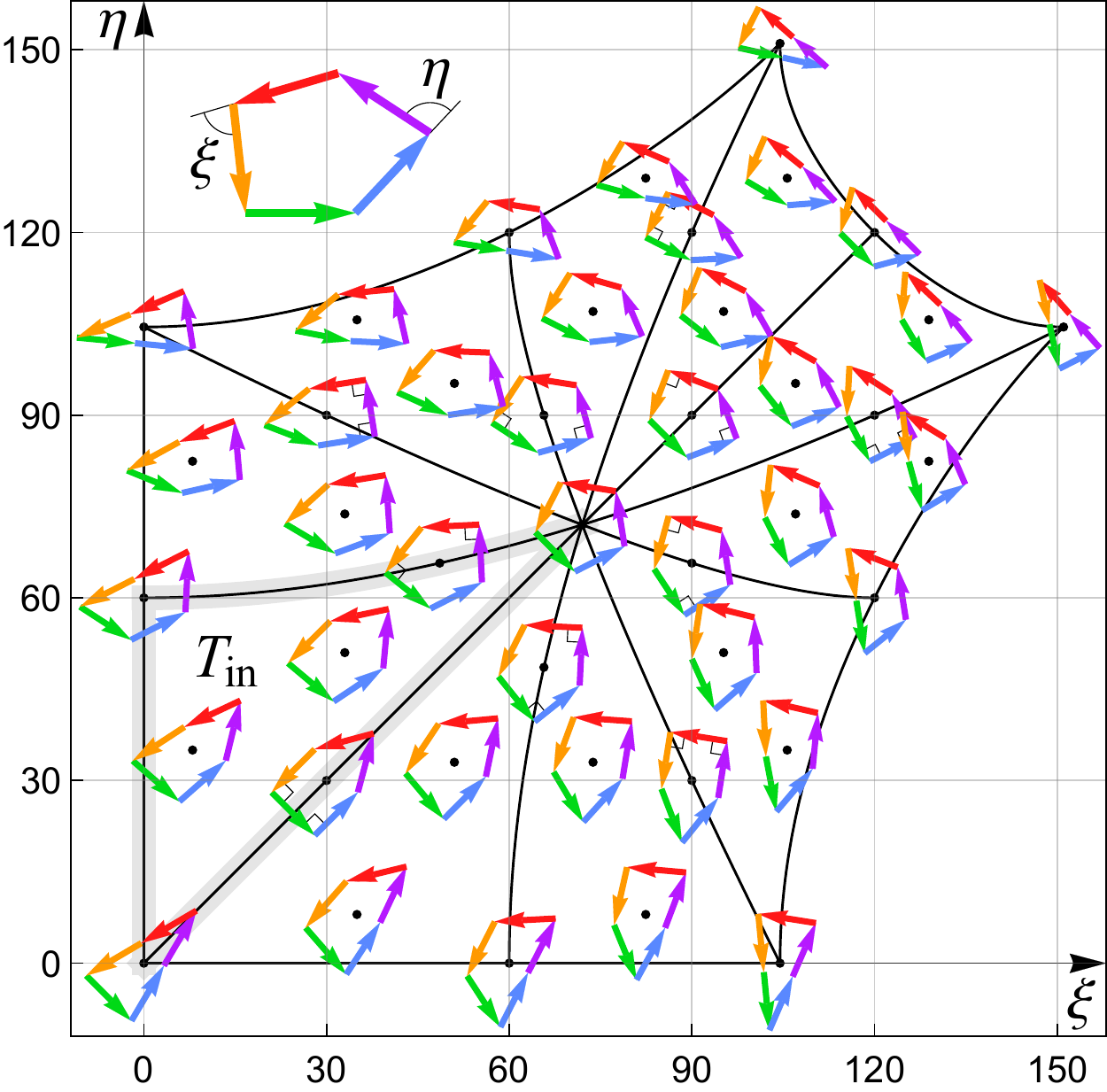}
	\end{center}
	\vspace*{-10 pt}
	\caption{The convex-and-\ccWise\ equilateral pentagons are here plotted using, as Cartesian coordinates $\xi$ and $\eta$, the external angles at the red-to-orange and blue-to-purple vertices.  The pentagons along the main diagonal are mirror symmetric on a green base, while those along each other curve through the center are mirror symmetric on a base of some other color.}
	\label{fig:GlimpseAng}
\end{figure}

\begin{figure}
	\begin{center}
		\qquad\includegraphics[scale=0.58]{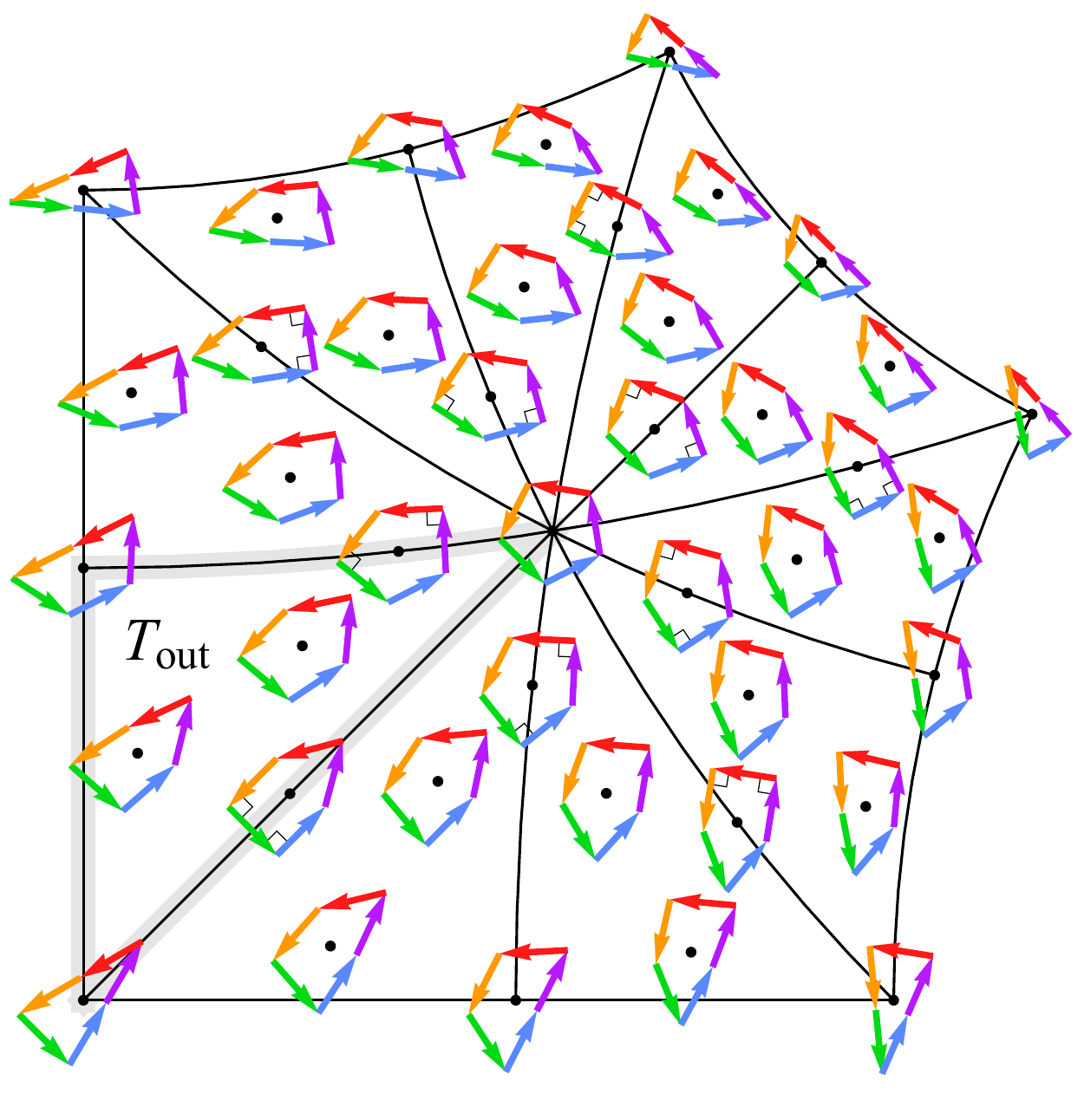}
	\end{center}
	\vspace*{-20 pt}
	\caption{The layout of Figure~\ref{fig:GlimpseAng} has here been adjusted to form a hyperbolic regular pentagon with right-angled vertices, drawn in a Poincar\'e projection.  These ten $(2,4,5)$\_triangles are $1/24^\text{th}$ of our conformal map of the polygon space $(E_5,e_5)$.}
	\label{fig:Glimpse}
\end{figure}

From among the equilateral pentagons in the Euclidean plane, Figure~\ref{fig:GlimpseAng} plots those that are convex and traversed counterclockwise (\ccWise), using two nonadjacent external angles as Cartesian coordinates.  That map breaks up into ten triangular regions based on which two of the five vertices have the smallest and largest angles.  (Those two vertices are always adjacent.)  In the more symmetric map of Figure~\ref{fig:Glimpse}, each of those ten regions has become a \emph{$(2,4,5)$\_triangle}, that is, a triangle in the hyperbolic plane whose vertex angles are $\frac{\pi}{2}$, $\frac{\pi}{4}$, and $\frac{\pi}{5}$.  

It\s easy to show the existence of the map in Figure~\ref{fig:Glimpse}.  The polygon space~$E_5$ of equilateral pentagons is an orientable, smooth $2$\_manifold of genus~4.  If we extend the five pentagon edges until all ten pairs of them intersect, the rates of change of those ten angles give us a Riemannian metric~$e_5$ on the manifold~$E_5$, a metric under which~$E_5$ has $240$ isometries.  When we equip~$E_5$ with the conformal structure implied by~$e_5$, we get a compact Riemann surface of genus~$4$ that has, as automorphisms, the $120$ isometries that preserve orientation.  But only one compact Riemann surface of genus $4$ has that many automorphisms:  the \emph{Bring sextic}.  The Bring sextic with its Poincar\'e metric embeds isometrically in the hyperbolic plane as a repeating pattern of $240$ $(2,4,5)$\_triangles, a pattern depicted in Figure~\ref{fig:Icosagon} on page~\pageref{fig:Icosagon}.  So there must be a conformal map of the manifold $(E_5,e_5)$ with that same structure, a map whose $240$ triangles are permuted transitively by the isometries of $(E_5,e_5)$.   Figure~\ref{fig:Glimpse} shows ten of that map\s $240$ triangles.

But those easy arguments take us only so far.  They tell us which equilateral pentagon lies at each vertex of a $(2,4,5)$\_triangle, which pentagons lie along each curve separating two triangles, and which lie interior to each triangle.  The labeled vertices act as signposts on a sketch of our entire map in Figure~\ref{fig:LabeledIcos} on page~\pageref{fig:LabeledIcos}.  But the full correspondence between pentagons and points takes more work.  For example, consider a \emph{child\s house}: a convex pentagon with two adjacent right angles.  Where on the main diagonal of Figure~\ref{fig:Glimpse} should the child\s house on a green base be plotted?  Answer: Near the hyperbolic line that forms the altitudes of the two adjacent map triangles, but below that line by a tiny distance --- about $0.48\%$ of a triangle\s hypotenuse, a distance just perceptible in Figure~\ref{fig:ScaleFactor} on page~\pageref{fig:ScaleFactor}.

The pentagon-to-point correspondence of Figure~\ref{fig:Glimpse} is a conformal map from the manifold $(E_5,e_5)$ to the Euclidean plane that takes each triangular region of~$E_5$ to the Poincar\'e projection of a $(2,4,5)$\_triangle.   We compute that map in stages.  We start with Figure~\ref{fig:GlimpseAng}, where the angles $\xi$ and $\eta$ at two nonadjacent vertices act as Cartesian coordinates.  The vertex angles of the ten regions in Figure~\ref{fig:GlimpseAng} differ from $\frac{\pi}{2}$, $\frac{\pi}{4}$, and $\frac{\pi}{5}$, showing that the coordinate system $(\xi,\eta)$ is not isothermal.   We find an initial pair of isothermal coordinates by following Gauss in converting the Beltrami partial differential equation over $\reals$ into an ordinary differential equation over $\complexes$.  That ODE is so complicated, unfortunately, that a solution in closed form seems hopeless; so we approximate a solution numerically.  This gives us a map of~$E_5$ in which the triangular regions have the proper vertex angles; but it takes a second conformal map, which we also approximate numerically, to convert one of those regions into the precise Poincar\'e projection of a $(2,4,5)$\_triangle, the other regions then falling into line automatically.

In this paper, after defining and defending the metric $e_5$, we deduce as much as we can from the easy arguments, postponing the hard work until the end.
\begin{itemize}
\item In Sections~\ref{sect:Strategy} and~\ref{sect:BringSextic}, we expand on this introduction, giving definitions, citing references, and recalling known results.  We equip each polygon space with a \emph{bag-of-edges} metric, the metric $e_5$ being the bag-of-edges metric for the polygon space $E_5$ of equilateral pentagons.
\item We defend the bag-of-edges metric in Section~\ref{sect:PhaseDifference} by showing that, under it, the loops formed by ``decomposable'' $n$\_gons are geodesics.  In Section~\ref{sect:Grassmannians}, we show that it accords with the Grassmannian metric that Hausmann and Knutson~\cite{Cantarella,HausmannKnutson} give to their larger moduli space of all planar $n$\_gons.
\item We determine in Sections~\ref{sect:VertexPentagons} and~\ref{sect:MapSketch} which equilateral pentagons lie along the boundaries of the triangular regions of~$E_5$ and which lie at their vertices, thereby exhausting what the easy arguments can tell us.
\item Section~\ref{sect:Dodeca} digresses by discussing a way to model $E_5$ as a surface in $3$\_space.
\item In Sections~\ref{sect:AngleCoords} through~\ref{sect:ConformalMap}, we buckle down to find the full pentagon-to-point correspondence.  This involves both symbolic and numeric computations that we relegate to the Mathematica appendix \texttt{EquilateralPentagonMap.nb} (or \texttt{.pdf}), whose $k^\text{th}$ section we herein cite as [A$k$].  To assist cartographers of pentagons with access to Mathematica, a function in that appendix [A9] converts the Cartesian coordinates $(\xi,\eta)$ of an equilateral pentagon in the triangle~$T_\text{in}$ of Figure~\ref{fig:GlimpseAng} into the Poincar\'e coordinates of that pentagon in~$T_\text{out}$ of Figure~\ref{fig:Glimpse} --- with any requested precision, up to $100$ digits.
\end{itemize}

\section{Our work plan and related work}
\label{sect:Strategy}

\subsection{Polygon spaces}

Given positive reals $\mathbf{r}=(r_1,\ldots,r_n)$, consider those $n$\_gons in the complex plane whose $k^\text{th}$ edge $\edge_k=r_k e^{\imagunit\theta_k}$ has length~$r_k$ and \emph{phase} (aka \emph{argument})~$\theta_k$, that phase lying in $S^1=\reals/2\pi\ints$.  For the $n$\_gon to be closed, its phases $(\theta_1,\ldots,\theta_n)$ must satisfy $\sum_kr_ke^{\imagunit\theta_k}=0$.  The \emph{polygon space} for the length vector $\mathbf{r}$ is the subset~$P_\mathbf{r}$ of $(S^1)^n/SO(2)$ satisfying that \emph{$\mathbf{r}$\_closure constraint}, where we mod out by the circle group $SO(2)$ to ignore rigid rotations of the entire $n$\_gon~\cite{Farber,Schutz}.  %(The polygon space $P_\mathbf{r}$ is empty when, for some $k$, we have $r_k>\sum_{i\ne k}r_i$.)

The $n$\_gons in a polygon space $P_\mathbf{r}$ have their vertices or edges labeled, so we know which edge is which, even when the length vector $\mathbf{r}$ is periodic.  Polygon spaces ignore translating and rotating their polygons.  Modding out by all of~$O(2)$, rather than just by~$SO(2)$, would ignore reflecting them as well.  The moduli space~$P_\mathbf{r}^{\scriptscriptstyle\updownarrow}$ that would result is half the size of $P_\mathbf{r}$, but is often nonorientable.

At a typical point in the polygon space $P_\mathbf{r}$, the tangent space to $(S^1)^n$ splits into three factors:  the tangent space to $P_\mathbf{r}$ at that point, a line of possible rigid rotations, and a plane of possible failure-to-close vectors.  The polygon space $P_\mathbf{r}$ is thus a compact, orientable, smooth $(n-3)$\_manifold, except possibly for finitely many singular points.  A singular point arises for each way of partitioning the lengths~$r_1$ through~$r_n$ into two subsets with the same sum, say subsets of cardinalities~$c$ and~$n-c$, since that enables the $n$\_gon to lie entirely along a single line.  At such a point, we get only a line of failure-to-close vectors, rather than a plane; so the local dimension of the polygon space jumps from $n-3$ to $n-2$.  The polygon space near that singular point is analytically isomorphic~\cite{KapMill95}\footnote{Kapovich and Millson only sketch the structure of singularities in~\cite{KapMill95}.  In~\cite{KapMill97,KapMill99}, however, they show that the singularities in a moduli space of spherical linkages are Morse using methods that apply to any space of constant curvature, including the plane.} to the origin in the zero set of a quadratic form on~$\mathbb{R}^{n-2}$ of signature $(c-1,n-c-1)$. 

Some other results: The homology groups of the polygon spaces are free abelian and of known rank~\cite[p.~16]{Farber}. Their cohomology rings are a bit subtle, even in the equilateral case~\cite{KamiyamaRing}.  But that ring structure almost determines $\mathbf{r}$, as Walker conjectured: Assuming that $\mathbf{r}$ has been sorted, that ring structure determines how the sum of each subset of $\mathbf{r}$ compares to the sum of its complement~\cite{Schutz}.

While all of this paper\s polygons are planar, the $n$\_gons in $\reals^d$ whose edge lengths are given by $\mathbf{r}$ form a polygon space also for $d>2$: the subset of $(S^{d-1})^n/SO(d)$ where the analogous $\mathbf{r}$\_closure constraint holds.  The case $d=3$ is especially pretty because the $2$\_sphere $S^2$ admits a symplectic form~\cite{HausmannKnutsonRings,KapMillSymp}.

\subsection{The bag-of-edges metric}

A tangent vector to a polygon space at some $n$\_gon is a recipe for varying the vertex angles of that $n$\_gon as a function of time.  One way to measure the length of such a vector would find the rates of change of the~$n$ vertex angles and take their $2$\_norm.  For greater symmetry, we extend the $n$ edges until all $\genfrac(){0 pt}{1}{n}{\raisebox{0.8pt}{$\scriptstyle 2$}}$ pairs of them intersect, finding the rates at which all of those angles change.  This makes our metric independent of the \emph{assembly order}, the linear order in which the $n$ edge vectors appear, tip to tail, going around the $n$\_gon.

Long edges matter more than short ones, so we weight the squared rate of change of an intersection angle by the lengths of the two intersecting edges:
\begin{equation}
\label{eq:BagOfEdges}
dp_\mathbf{r}^2=\sum_{i<\mkern 1mu j}r_i r\subj (d\theta\subj-d\theta_i)^2.
\end{equation}
We equip the polygon space $P_\mathbf{r}$ with the \emph{bag-of-edges metric}, the Riemannian metric $p_\mathbf{r}$ given by this first fundamental form.

The metric $p_\mathbf{r}$ behaves well on $n$\_gons that are ``decomposable''\kern -1pt.  The~$n$ directed edges of any $n$\_gon always sum to zero, but it may happen that some proper subset of them also sums to zero, in which case the complementary subset must sum to zero as well.  We call such an $n$\_gon \emph{decomposble}.  Any decomposable $n$\_gon lies on a loop that is formed by rotating the edges of the first subset in lockstep while keeping the complementary edges fixed.  As we show in Section~\ref{sect:PhaseDifference}, our bag-of-edges metric $p_\mathbf{r}$ has the nice property that any such \emph{loop of decomposables} is a geodesic of the Riemannian manifold $(P_\mathbf{r},p_\mathbf{r})$.

A stronger argument in favor of the metric $p_\mathbf{r}$ comes from the Grassmannian approach of Hausmann and Knutson~\cite{Cantarella,HausmannKnutson}.  Suppose that we allow the edge lengths of our $n$\_gons to vary along with their vertex angles; we fix only their perimeter, so as to ignore magnification.  The moduli space of all $n$\_gons with some fixed perimeter has dimension $2n-4$.  Hausmann and Knutson model that space as the Grassmannian~$G(2,n)$ of $2$\_dimensional linear subspaces of $\reals^n\mkern -1mu$.  This puts a metric on that moduli space --- a metric $g_n$ that, like our $p_\mathbf{r}$, ignores the assembly order of the edges.  For any length vector~$\mathbf{r}$ that sums to the fixed perimeter, the polygon space~$P_\mathbf{r}$ sits inside that Grassmannian, roughly speaking.  We discuss further in Section~\ref{sect:Grassmannians}, showing that our bag-of-edges metric~$p_\mathbf{r}$ differs from the metric that~$P_\mathbf{r}$ inherits as a subset of that Grassmannian just by a constant factor: $dp_\mathbf{r}^2=8\,dg_n^2$.

\subsection{The equilateral case}

The length vector $\mathbf{e}(n)=\bigl(\frac{1}{n},\ldots,\frac{1}{n}\bigr)$ gives us the polygon space $P_{\mathbf{e}(n)}$ of equilateral $n$\_gons with unit perimeter.  For brevity, we denote that polygon space $P_{\mathbf{e}(n)}$ simply as $E_n$ and its bag-of-edges metric $p_{\mathbf{e}(n)}$ as $e_n$.  

Let the symmetric group $S_n$ permute the assembly order of the edges, while the cyclic group $C_2$ optionally reflects the entire $n$\_gon.  All $2n!$ of the actions of $S_n\times C_2$ on~$E_n$ are then isometries of the Riemannian manifold $(E_n, e_n)$.

Are those isometries distinct?  Once $n\ge 5$, there exist labeled equilateral $n$\_gons that are fixed, under the actions of $S_n\times C_2$, only by the action arising from the identity; so the $2n!$ isometries are then distinct.  They are not distinct for $n\le 3$, but they turn out to be also for $n=4$ --- an interesting example.

\begin{figure}
	\begin{center}
	\includegraphics[scale=0.6]{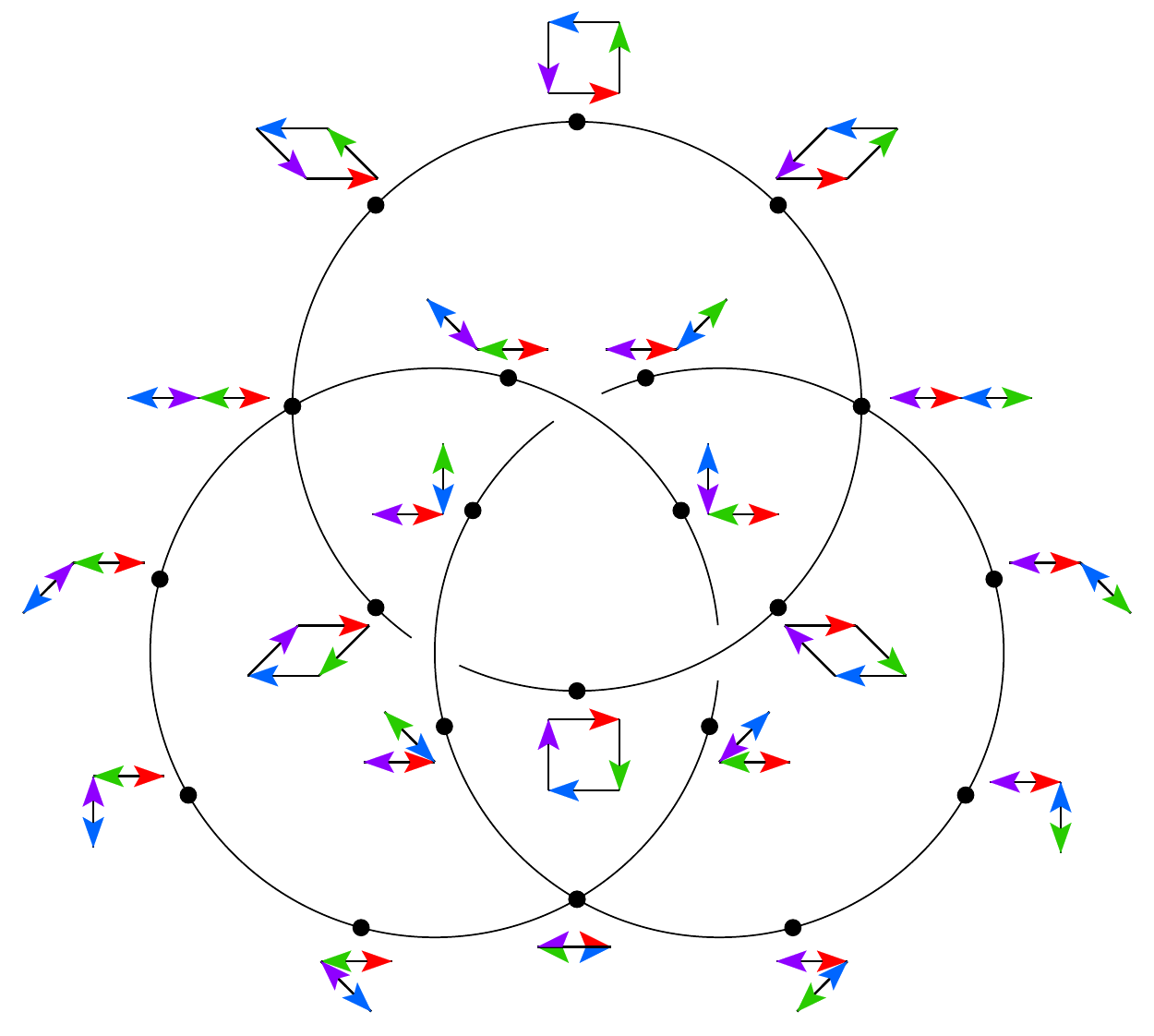}
	\end{center}
	\vspace{-15pt}
	\caption{The polygon space $E_4$ of equilateral quadrilaterals consists of three loops:  one of rhombs and two of back-and-forth L\s.  Each pair of loops intersects once (orthogonally under $e_4$ [A1], though not in this sketch) at a quadrilateral that lies entirely along a line, those intersections being the three singular points of $E_4$.}
	\label{fig:Tetravecs}
\end{figure}

Figure~\ref{fig:Tetravecs} depicts the polygon space $E_4$ of equilateral quadrilaterals.  While lots of quadrilaterals don\t decompose, every equilateral quadrilateral decomposes into two digons.  The space $E_4$ consists of three loops of such decomposables, loops that are geodesics under $e_4$.  It has three singular points, one for each way of partitioning the four edges into two pairs.  Each looks locally like the origin in the zero set of a quadratic form on $\reals^2$ of signature~$(1,1)$, a form such as $y^2-x^2$.  The~$2\cdot 4!=48$ isometries of the manifold $(E_4,e_4)$ combine permuting the three loops (in $6$ ways) with swapping the halves of some subset of them (in $8$ ways).\footnote{The isometries of $(E_4,e_4)$ illustrate the isomorphism between the product group $S_4\times C_2$ and the hyperoctohedral group $S_2\wr S_3=(C_2)^3\rtimes S_3$, where the $S_3$ in that semidirect product permutes the three loops while each $C_2$ is optionally swapping the halves of one of them.}

\subsection{Equilateral pentagons}

Our focus is the polygon space $E_5$ of equilateral pentagons.  Topologically, $E_5$ has been shown --- frequently~\cite{CurtisSteiner,Havel,KamiyamaNoMorse,KapMill95,ShimamotoVanderwaart} --- to be a compact, orientable, smooth $2$\_manifold of genus~$4$.  Klaus and Kojima~\cite{KlausKojima} study $E_5$ as an algebraic variety.   In this paper, we study $E_5$ using differential geometry.

Tong and Dullin~\cite{TongDullin} also studied $E_5$ geometrically, their Figure~5~\cite[p.~975]{TongDullin} being a nonhyperbolic analog of our Figure~\ref{fig:LabeledIcos}.  They consider planar equilateral pentagons with unit-length massless rods as their edges and with unit-mass rotary actuators as their vertices.  They study how such a linkage can achieve a global rotation by applying an appropriate pattern of actuator torques, thus giving some insight into how a falling cat manages to reorient itself so as to land on its feet.  But their variant of the polygon space $E_5$ is only $1/12^\text{th}$ as symmetric as ours:  Of the $120$ permutations of the assembly order of the actuators, only the $10$ dihedral ones preserve the dynamics of their linkage.

\subsection{The Bring sextic}

Equipping $E_5$ with our bag-of-edges metric $e_5$ gives a Riemannian $2$\_manifold $(E_5,e_5)$ with $2\cdot 5!=240$ isometries.  Suppose that we forget how the metric~$e_5$ is scaled, remembering only the conformal structure $e_5^\smangle$ that it gives to $E_5$, and that we orient~$E_5$ one way or the other, denoting the result~$E_5^+$.  (Section~\ref{subsect:Orientation} proposes a preferred orientation, but either is fine for now.)  We get a Riemann surface $\RiemannEfive=(E_5^+,e_5^\smangle)$ that is compact, has genus $4$, and has at least $120$ automorphisms:  the isometries of $(E_5,e_5)$ that preserve orientation.

Among the compact Riemann surfaces of genus~$3$, the famous \emph{Klein quartic} has the most automorphisms, at $168$.  Among those of genus~$4$, the maximum number is $120$~\cite{RodriguezGonzalez,Kimura}, and the \emph{Bring sextic}~\cite[$\Phi_{50}$ on p.~32]{Kimura} is the unique surface, up to conformal equivalence, that has~$120$ automorphisms.  (We say ``the Bring sextic'' by analogy with the standard phrase ``the Klein quartic''\kern -1pt, thus obviating the need to choose between ``the Bring curve'' over $\complexes$ and ``the Bring surface'' over $\reals$.  Back in genus~$2$, this analogy suggests saying ``the Bolza quintic''\kern -1pt.)  

The Bring sextic, which we denote $\Bsextic$, can be constructed as the points in complex projective $4$\_space whose five homogeneous coordinates, their squares, and their cubes all sum to zero.  The first of those constraints cuts out a projective $3$\_space, while the others define a quadric and a cubic surface in that $3$\_space.  Their complete intersection is a nonsingular sextic space curve that Klein named after the Swedish mathematician Erland Samuel Bring (1736--1798).   Its $120$ automorphisms are the permutations of the five homogeneous coordinates. 

By the uniformization theorem, any compact Riemann surface of genus at least~$2$ has a unique \emph{Poincar\'e metric} that is consistent with its conformal structure and makes its Gaussian curvature constant at~$-1$.  Under that Poincar\'e metric, the surface embeds isometrically as a repeating pattern in the hyperbolic plane.  As we review in Section~\ref{sect:BringSextic} and depict in Figure~\ref{fig:Icosagon}, that repeating pattern for the Bring sextic $\Bsextic$ has, as one of its fundamental domains, a convex hyperbolic icosagon that is tiled by $240$ $(2,4,5)$\_triangles.  

Given the genus and symmetries of our Riemann surface $\RiemannEfive=(E_5^+,e_5^\smangle)$, it must be conformally equivalent to the Bring sextic $\Bsextic$.  We thus deduce at once that our conformal map of $(E_5,e_5)$ will have the structure of that repeating icosagon, with the $240$ isometries of $(E_5,e_5)$ transitively permuting the triangles.

\subsection{Labeling the vertices}

We refer to the $(2,4,5)$\_triangles that tile the Bring sextic~$\Bsextic$ as \emph{tracts}.  Since the surface $\RiemannEfive=(E_5^+,e_5^\smangle)$ shares that structure, $E_5$ is also partitioned into tracts.

The curves that bound the tracts on $\Bsextic$ are loops that are pointwise fixed by self-inverse antiautomorphisms.  Similarly, the curves that bound the tracts on $(E_5,e_5)$ are loops pointwise fixed by self-inverse, orientation-reversing isometries.  In Section~\ref{sect:MapSketch}, we determine those isometries and the equilateral pentagons that they fix.  Since any pentagon that lies at a tract vertex is fixed by several such isometries, we can determine precisely which pentagons lie at the $24$ $\frac{\pi}{5}$\_vertices of~$E_5$, at its $30$ $\frac{\pi}{4}$\_vertices, and at its $60$ $\frac{\pi}{2}$\_vertices.   We thereby sketch our entire conformal map in Figure~\ref{fig:LabeledIcos}:  a drawing of the repeating icosagon in which the tract vertices are labeled with their equilateral pentagons.

The bad news is that we don\t yet know the pentagon-to-point correspondence either within a tract or along the boundary curve separating two tracts. The easy arguments don\t suffice for that.

Before tackling that challenge, however, we pause to describe an attractive way to model~$E_5$ as a surface in Euclidean $3$\_space, since that model is a big help in clarifying how~$E_5$ is connected up globally.   Weber~\cite{Weber} suggests modeling the Bring sextic using the self-intersecting polyhedron called a \emph{dodecadodecahedron}, and Section~\ref{sect:Dodeca} points out that this works particularly well for $E_5$.

\subsection{The pentagon-to-point correspondence}

We then begin the hard work of computing the precise pentagon-to-point correspondence of our conformal map.  In Section~\ref{sect:AngleCoords}, we adopt the vertex angles~$\xi$ and~$\eta$ from Figure~\ref{fig:GlimpseAng} as our coordinates, computing both the first fundamental form~$de_5^2$ of the bag-of-edges metric $e_5$ and the Gaussian curvature~$K\mkern -1mu$ of the manifold $(E_5,e_5)$ in those coordinates.  We find that~$K$ varies over $[-11/3\dd-1]$.

Any planar map of a Riemannian $2$\_manifold is the coordinate plane of some coordinate system.  When the map is conformal, that coordinate system is said to be \emph{isothermal} (a term from heat diffusion).  Isothermal systems are the solutions of the \emph{Beltrami} partial differential equation.  Following Gauss, we solve the Beltrami by converting it into an ordinary differential equation over~$\complexes$ in Section~\ref{sect:IsoCoords}.  Selecting the second of two candidate solutions, we get, as a conformal image of a tract of $(E_5,e_5)$, a triangular region $T_\text{cnf}$ in~$\complexes$ with a straight side along the real axis and with two curved sides that we know only numerically.  

Section~\ref{sect:ConformalMap} conformally maps $T_\text{cnf}$ to a Poincar\'e projection of a $(2,4,5)$\_triangle.  A small change of shape suffices, so we approximate that conformal map using a polynomial of one complex variable with real coefficients.  This gives us, at last, the full pentagon-to-point correspondence of our conformal map.

The Gaussian curvature of the manifold $(E_5,e_5)$ is not constant, but standard theory tells us that the metric $e_5$ can be rescaled by a unique smooth function $s\colon E_5\to\mathbb{R}_{>0}$ so that the rescaled manifold $(E_5,s^2e_5)$ has Gaussian curvature~$-1$ everywhere.  We approximate that rescaling $s$ by taking our pentagon-to-point correspondence and backing out the effect of its Poincar\'e projection.  Perhaps surprisingly, $s$ turns out to vary only modestly, with $\max(s)/\min(s)\approx 1.04$ --- so a gentle rescaling suffices to normalize the curvature of $E_5$.

\subsection{The equiangular analog}

In the theory of polygon spaces, the lengths of the edges are held fixed while the vertex angles are allowed to vary.  There is an intriguingly analogous theory in which the vertex angles are held fixed while the edge lengths vary.  That analogous theory doesn\t arise as often in applications, and it won\t arise elsewhere in this paper; but it has the advantage of being technically simpler.  We discuss it here to clarify the relationships between the two theories.

Polygon spaces are curved because fixing the edge lengths of an $n$\_gon constrains its vertex angles nonlinearly; but the constraints that fixed vertex angles impose on the edge lengths are linear.  Given positive numbers $\bm{\rho}=(\rho_1,\ldots,\rho_n)$ that sum to~$1$, consider those convex $n$\_gons in which the $k^\text{th}$ vertex has external angle~$2\pi\rho_k$.  When we scale those $n$\_gons to have a fixed perimeter, they form a Euclidean polytope of dimension~$n-3$, each facet of which consists of those $n$\_gons in which a particular edge has shrunk to length zero~\cite{RamshawSaxe}.  We could view those fixed-perimeter $n$\_gons as a subset of the Grassmannian $G(2,n)$~\cite{Cantarella,HausmannKnutson}.  But the prettiest option scales those $n$\_gons to have unit area, as both Bavard and Ghys~\cite{BavardGhys} and Kojima and Yamashita~\cite{KojYam} suggest.  When scaled that way, those $n$\_gons form a hyperbolic polytope with lots of right angles, called a \emph{truncated orthoscheme}~\cite{Fillastre}.
(Thurston~\cite[pp.~72--73]{Thurston} assigns that hyperbolic option as Problem 2.3.12.)

The cyclic order of the $n$ fixed angles is itself fixed in such a Bavard--Ghys polytope.  Kapovich and Millson~\cite{KapMill95} allow $n$ vertices with the external angles in~$\bm{\rho}$ to occur in any order, going around the $n$\_gon.  The moduli space that they get is a hyperbolic cone-manifold~$K_{\bm{\rho}}$ tiled by $(n-1)!$ Bavard--Ghys polytopes~\cite{Fillastre}.  Using Schwarz--Christoffel maps, they define a homeomorphism from the fixed-lengths polygon space $P_{\bm{\rho}}$ to their fixed-angles cone-manifold $K_{\bm{\rho}}$.  The parameters~$\rho_k$ here specify the fixed edge lengths in $P_{\bm{\rho}}$, but the fixed vertex angles in $K_{\bm{\rho}}$.

We know about the polygon space $P_{\mathbf{e}(5)}=E_5$ of equilateral pentagons.  The cone-manifold $K_{\mathbf{e}(5)}$ of equiangular pentagons turns out to be tiled by $24$ regular, right-angled hyperbolic pentagons, with four around each vertex.  The polygon space $P_{\mathbf{e}(5)}$ has that same structure, with ten of its $240$ tracts forming each of the hyperbolic pentagons.  Thus, if the Kapovich--Millson homeomorphism from $P_{\mathbf{e}(5)}$ to $K_{\mathbf{e}(5)}$ became conformal when~$P_{\mathbf{e}(5)}$ was given the bag-of-edges metric~$p_{\mathbf{e}(5)}=e_5$, our hard work in Sections~\ref{sect:AngleCoords}--\ref{sect:ConformalMap} would be unnecessary.  Far from being conformal, however, their homeomorphism isn\t even a diffeomorphism.  Consider an equilateral pentagon in~$P_{\mathbf{e}(5)}$ that, for simplicity, is convex.  As one of its external angles~$\rho$ goes to zero, the corresponding side length~$r$ of the equiangular pentagon in $K_{\mathbf{e}(5)}$ also goes to zero; but~$r$ goes to zero, not like $\rho$, but like $\sqrt[5]{\rho}$.

\section{The repeating icosagon}
\label{sect:BringSextic}

\begin{figure}
	\begin{center}
		\includegraphics[scale=0.8]{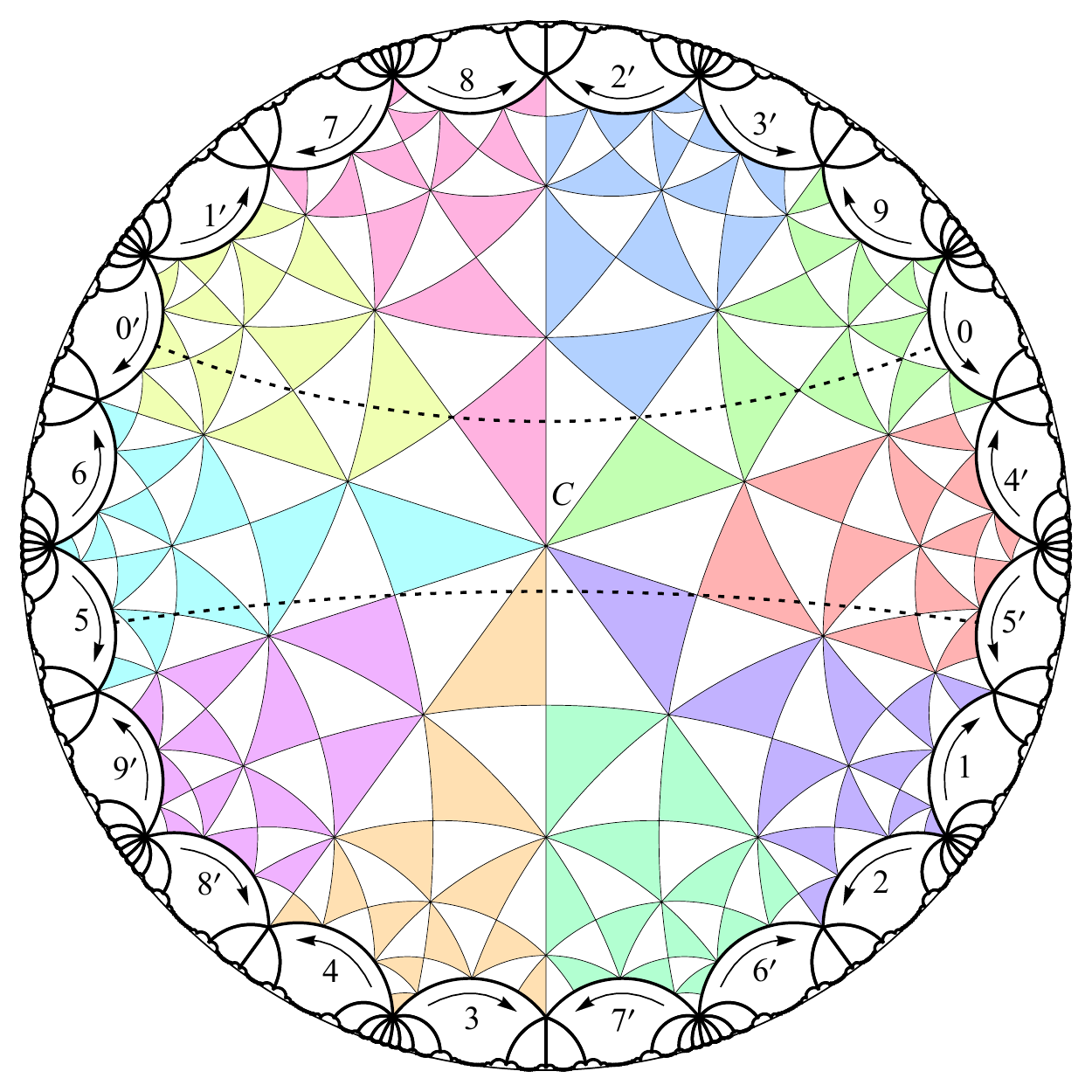}
		\begin{picture}(0,0)
     	\end{picture}
	\end{center}
	\vspace*{-15 pt}
	\caption{This icosagon, consisting of $240$ $(2,4,5)$\_triangles, is a fundamental polygon for the repetitions of the hyperbolic embedding of the Bring sextic $\Bsextic$.  The translates of the icosagon that touch it are outlined in black. (Its ten merged pairs of edges are here numbered in \ccWise\ order as they leave its merged $\frac{\pi}{5}$ vertex.)}
	\label{fig:Icosagon}
\end{figure}

To complete the relevant background, let\s recall how the Bring sextic $\Bsextic$, equipped with its Poincar\'e metric, embeds isometrically in the hyperbolic plane.

Twenty-four $(2,4,5)$\_triangles can reflect across their edges to form a hyperbolic regular quadrilateral with vertex angles of $\frac{\pi}{5}$; we call such a quadrilateral a \emph{precinct}.  Figure~\ref{fig:Icosagon} shows ten precincts, differently tinted, surrounding the point $C$ and filling out an equilateral icosagon whose vertex angles alternate between $\frac{\pi}{5}$ and $\frac{2\pi}{5}$.  

For $0\le k\le 9$, there is a hyperbolic translation that carries the edge~$k$ in Figure~\ref{fig:Icosagon} to the edge~$k'$.  The upper dotted curve is the axis of the translation that takes $0$ to~$0'$, while the lower takes $5$ to $5'$; the other even and odd $k$ are symmetric.  Those translations generate a group of orientation-preserving isometries of the hyperbolic plane, and the images of the icosagon under that group tile the plane.  Identifying all of the translates merges the ten $\frac{\pi}{5}$\_vertices, the five $\frac{2\pi}{5}$\_vertices at the heads of the even arrows, and the five at the~heads of the odd arrows.  After those merges, the ten precincts share their four vertices.

The icosagon in Figure~\ref{fig:Icosagon}, which is the Dirichlet polygon of $C\mkern -1mu$, is a fundamental polygon for the embedding of the Bring sextic $\Bsextic$ as a repeating pattern in the hyperbolic plane~\cite{RieraRodriguez,Weber}.  Since sewing the ten precincts together gives a surface with $4$ vertices, $20$ edges, and $10$ faces, its Euler characteristic $V-E+F$ is $-6$, which is correct for a surface of genus $4$.

As we have discussed, that icosagon must be a fundamental polygon also for the repetitions of our conformal embedding of the polygon space~$(E_5,e_5)$.

\paragraph{Tracts}

We call each $(2,4,5)$\_triangle in either $\Bsextic$ or $E_5$ a \emph{tract} --- either an \emph{\Ltract} or an \emph{\Rtract} according as the turn from its long leg to its short leg is to the left (as in the letter L) or to the right.  In Figure~\ref{fig:Icosagon}, the \Rtract{}s have the tints.

A \emph{decatract} is the ten tracts that surround a $\frac{\pi}{5}$\_vertex, such as $C\mkern -1mu$.  Each is a hyperbolic regular, right-angled pentagon bounded by five pairs of short legs.  Since $\Bsextic$ has $240$ tracts overall, it is tiled by~$24$ decatracts.  

\paragraph{Geodesic loops}

Each curve on $\Bsextic$ that separates adjacent tracts is the geodesic loop that is pointwise fixed by some self-inverse antiautomorphism.\footnote{Self-inverse antiautomorphisms are so important in the study of Riemann surfaces that the word ``symmetry'' is often restricted to refer only to them --- for example, in the title of~\cite{strictSymmetries}.}   Of the $120$ antiautomorphisms of $\Bsextic$, there are $26$ that are self-inverse:
\begin{itemize}
\item Just conjugating all five homogeneous coordinates gives a self-inverse antiautomorphism; but it has no fixed points, so no geodesic loop results.
\item Swapping one pair of homogeneous coordinates and conjugating all of them fixes one of the $10$ \emph{short-leg loops}, each consisting of $12$ short legs.
\item Swapping two disjoint pairs of coordinates and conjugating fixes one of the $15$ \emph{hypot-long-leg loops}, each consisting of $8$ hypotenuses and $8$ long legs.
\end{itemize}

Among other geodesic loops, the axis of the translation that takes edge $0$ to $0'$ is one of $20$ \emph{altitude loops}, each of which consists of $12$ tract altitudes.  The line segment joining the midpoints of nonadjacent sides of a decatract is the summit of a Saccheri quadrilateral; and the axis of the $5$\_to\_$5'$ translation is one of $30$ \emph{Saccheri-summit loops}, each of which consists of four such summits.\footnote{The lengths of the altitude, Saccheri-summit, and short-leg loops, which are about $4.603$, $4.795$, and $6.368$, start off the length spectrum of $\Bsextic$.  The length of the hypot-long-leg loops, about $11.755$, is dozens of places farther out in that spectrum.  

Removing a single loop of any of those four types from $\Bsextic$ does not disconnect it.  But there are sets of six altitude loops and three Saccheri-summit loops that are disjoint; and removing nine such loops from~$\Bsextic$ is an elegant way to cut it up into six isometric pairs of pants, with waists just a bit longer than their cuffs.  The graph encoding how the six pants are sewn together is~$K_{3,3}$, the waist twists are $1/4$ (side seams aligned with crotch), and the cuff twists are $1/6$.}

%\footnote{The lengths of the altitude, Saccheri, and short-leg loops are $2\arccosh\bigl((9+5\smash{\sqrt{5}})/4\bigr)\approx 4.603$, $2\arccosh\bigl((11+5\smash{\sqrt{5}})/4\bigr)\approx 4.795$, and $2\arccosh\bigl((13+5\smash{\sqrt{5}})/2\bigr)\approx 6.368$, and those three values start off the length spectrum of $\Bsextic$.  The hypot-long-leg loops have length $2\arccosh(89+40\smash{\sqrt{5}})\approx 11.755$, which is dozens of places farther out in that spectrum.}

\section{Loops of decomposables}
\label{sect:PhaseDifference}  

We are equipping each polygon space $P_\mathbf{r}$ with our \emph{bag-of-edges} metric $p_\mathbf{r}$.  One argument in favor of that choice of metric is that any loop of decomposables is a geodesic of the Riemannian manifold $(P_\mathbf{r}, p_\mathbf{r})$, as we here show.

The polygon space $P_\mathbf{r}$ is defined to be the subset of the quotient $(S^1)^n/SO(2)$ in which $(\theta_1,\ldots,\theta_n)$ satisfies the closure constraint $\sum_k r_k e^{\imagunit \theta_k}=0$, where we mod out by $SO(2)$ to ignore rigid rotations.  By working with the $\binom{n}{2}$ phase differences, rather than with the $n$ phases directly, we can equivalently define $P_\mathbf{r}$ to be a subset of a larger torus without needing to mod out by $SO(2)$.

We define the \emph{all-pairs torus} to be the Cartesian product $\smash{(S^1)^{\binom{n}{2}}}$ of $\binom{n}{2}$ copies of the circle $S^1$.  The \emph{phase-difference model} of the polygon space $P_\mathbf{r}$, which we denote $P_\mathbf{r}^\Delta$, is the subset of the all-pairs torus containing the points whose $(i,j)^\text{th}$ coordinate $c_{i,j}$, for $1\le i<j\le n$, is the difference $c_{i,j}=\theta\subj-\theta_i$ for some phases $(\theta_1,\ldots,\theta_n)$ that satisfy the closure constraint.   The model $P_\mathbf{r}^\Delta$ lies in the flat subtorus of the all-pairs torus where $c_{i,j}=c_{i,i+1}+\cdots+c_{j-1,j}$.  That subtorus is of dimension $n-1$; and it is affinely isomorphic to $(S^1)^n/SO(2)$ by an isomorphism that takes the phase-difference model $P_\mathbf{r}^\Delta$ to the standard model $P_\mathbf{r}$. 

We exploit the phase-difference model by giving the all-pairs torus the cubical metric of Euclidean $\binom{n}{2}$\_space, but with its $(i,j)^\text{th}$ coordinate axis scaled by $\sqrt{\mystrut{3 pt}\smash{r_ir\subj}}$, thus viewing the all-pairs torus as a cuboid with opposite facets identified.  The bag-of-edges metric $dp_\mathbf{r}^2=\sum_{i<\mkern 1mu j}r_i r\subj (d\theta\subj-d\theta_i)^2$ 
is then precisely the metric that the polygon space $P_\mathbf{r}^\Delta$ inherits as a submanifold of that cuboidal torus.

Now, consider an $n$\_gon that decomposes into a $k$\_gon and an $(n-k)$\_gon, for some $0<k<n$.  As we rotate the $k$ edges in lockstep, each phase difference $\theta\subj-\theta_i$ varies at a constant rate of either $+1$,~$0$, or~$-1$.  The loop of decomposables is thus a straight line through the all-pairs torus that remains entirely inside the polygon space $P_\mathbf{r}^\Delta$.  That straight line is a geodesic of the all-pairs torus and hence also a geodesic of the Riemannian manifold $(P_\mathbf{r},p_\mathbf{r})$.

\section{The Grassmannian moduli space for \texorpdfstring{$\bm{n}$}{n}\_gons}
\label{sect:Grassmannians}  

A weightier argument in favor of our bag-of-edges metric $p_\mathbf{r}$ is that it accords with the Grassmannian approach of Hausmann and Knutson~\cite{HausmannKnutson}, which is nicely explained by Cantarella, Needham, Shonkwiler, and Stewart~\cite{Cantarella}.  They fix only the perimeter of their $n$\_gons, allowing the individual edge lengths to vary along with the vertex angles.  The resulting moduli space has dimension $2n-4$, and they model it using the Grassmannian $G(2,n)$ of $2$\_dimensional linear subspaces of~$\reals^n\mkern -1mu$.

Some fine points: To simplify their formulas, they choose their fixed perimeter to be~$2$.  When they use points in $G(2,n)$ to model $n$\_gons of perimeter~$2$, those $n$\_gons are identified under reflection, as well as under translation and rotation.  To avoid ignoring reflection, they use the oriented Grassmannian~$\widetilde G(2,n)$.  Finally, there is a multiple-covering issue.  For a length vector~$\mathbf{r}$ that sums to $2$, the subset of $\widetilde G(2,n)$ that models $n$\_gons with edges of length $\mathbf{r}$ is, not the polygon space $P_\mathbf{r}$ itself, but a $2^{n-1}$\_fold cover\footnote{A $2^{n-1}$\_fold cover even when $P_\mathbf{r}$ has singular points; and regular, with $(C_2)^n/C_2\cong (C_2)^{n-1}$ as its group of covering transformations~\cite[p.~180]{HausmannKnutson}.  If we ignore reflection, though, the singular points of the half-sized $P_\mathbf{r}^{\scriptscriptstyle\updownarrow}$ must typically be deleted for $G(2,n)$ to provide a cover of what\s left.} of $P_\mathbf{r}$.  Despite this, it still makes sense to compare a metric on~$P_\mathbf{r}$ with the standard metric on $\widetilde G(2,n)$, which we denote $g_n$.

In this section, we show that $dp_\mathbf{r}^2=8\,dg_n^2$ --- with the proviso that~$p_\mathbf{r}$ is defined only for those changes to the phases of the $n$\_gon that preserve its closure to first order when its edge lengths are held fixed.  In contrast, $g_n$ would allow any changes to its phases and edge lengths that preserve its closure and perimeter to first order.  As for the factor of $8$, magnifying the $n$\_gons to raise their perimeter from $1$ to $2$ accounts for a factor of $4$; but the remaining factor of $2$ is more obscure.

(Readers not interested in Grassmannians may skip the rest of this section.)

Consider some $n$\_gon $\bm{\theta}=(\theta_1,\ldots,\theta_n)$ in a polygon space $P_\mathbf{r}$ where $\mathbf{r}$ sums to~$2$.  By the $\mathbf{r}$\_closure constraint, we have $\sum_k r_k\cos\theta_k=\sum_k r_k\sin\theta_k=0$.  The points in~$\widetilde G(2,n)$ that cover $\bm{\theta}$ are~\cite{Cantarella} the oriented $2$\_spaces $W\!=\Span(\mathbf{u},\mathbf{v})$, where
\[u_k=\sigma_k\sqrt{r_k}\cos(\theta_k/2)\qquad\text{and}\qquad v_k=\sigma_k\sqrt{r_k}\sin(\theta_k/2)\]
for $1\le k\le n$.  Each $\sigma_k$ is either $\pm 1$, and those choices determine which covering point $W$ we get.  There are only $2^{n-1}$ covering points because negating all $n$ of the~$\sigma_k$ gives a different basis for the same oriented $2$\_space $W\mkern -2mu$.  And all of those bases $(\mathbf{u},\mathbf{v})$ are orthonormal because of the trigonometric identities
\begin{align*}
\cos(\theta/2)^2&=(1+\cos\theta)/2\\
\sin(\theta/2)^2&=(1-\cos\theta)/2\\
\cos(\theta/2)\sin(\theta/2)&=(\sin\theta)/2.
\end{align*}
 
A tangent vector $d\bm{\theta}$ to $P_\mathbf{r}$ at the $n$\_gon $\bm{\theta}$ is determined by infinitesimal phase changes~$d\theta_k$ with $\sum_k r_k\cos(\theta_k+d\theta_k)=\sum_k r_k\sin(\theta_k+d\theta_k)=0$, from which it follows that $\sum_k r_k\cos(\theta_k)\,d\theta_k=\sum_k r_k\sin(\theta_k)\,d\theta_k=0$. 
The squared length of $d\bm{\theta}$ under the bag-of-edges metric $p_\mathbf{r}$ is, from~(\ref{eq:BagOfEdges}),
\begin{align}
\notag 
\sum_{i<\mkern 1mu j}r_ir\subj(d\theta\subj-d\theta_i)^2&=
\frac{1}{2}\sum_{i,\,j}r_ir\subj(d\theta_i^2-2\,d\theta_i\,d\theta\subj
+d\theta\subj^2)\\
\notag
&=\Bigl(\sum_ir_i\Bigr)\Bigl(\sum_j r\subj\, d\theta\subj^2\Bigr)-
\sum_{i,\,j}r_ir\subj\,d\theta_i\,d\theta\subj\\
\label{eq:TangVecBagLength}
&=2\sum_kr_k\,d\theta_k^2-\Bigl(\sum_k r_k\,d\theta_k\Bigr)^{\!2}.
\end{align}

Under the Grassmannian metric $g_n$, the tangent vector $d\bm{\theta}$ is a linear map from $W$ to $W^\perp$, and its squared length (squared Frobenius norm) is the sum of the squared lengths of the projections of $d\mathbf{u}$ and of $d\mathbf{v}$ into $W^\perp$, where
\[du_k=-{\textstyle\frac{1}{2}}\sigma_k\sqrt{r_k}\sin(\theta_k/2)\,d\theta_k\qquad\text{and}\qquad dv_k={\textstyle\frac{1}{2}}\sigma_k\sqrt{r_k}\cos(\theta_k/2)\,d\theta_k.\]
Since $d\mathbf{u}\cdot\mathbf{u}$ simplifies to zero, the projection of $d\mathbf{u}$ into $W^\perp$ is $d\mathbf{u}-(d\mathbf{u}\cdot\mathbf{v})\mathbf{v}$.  Dotting that with itself, its squared norm is 
$d\mathbf{u}\cdot d\mathbf{u}-2(d\mathbf{u}\cdot\mathbf{v})^2+(d\mathbf{u}\cdot\mathbf{v})^2(\mathbf{v}\cdot\mathbf{v})$, which simplifies to $d\mathbf{u}\cdot d\mathbf{u}-(d\mathbf{u}\cdot\mathbf{v})^2$.  
The squared norm of the projection of $d\mathbf{v}$ into $W^\perp$ is $d\mathbf{v}\cdot d\mathbf{v}-(d\mathbf{v}\cdot\mathbf{u})^2$ similarly, and we then calculate that
\[d\mathbf{u}\cdot d\mathbf{u}+d\mathbf{v}\cdot d\mathbf{v}=\frac{1}{4}\sum_kr_k\,d\theta_k^2,
\quad\text{while}\quad d\mathbf{v}\cdot\mathbf{u}=-d\mathbf{u}\cdot\mathbf{v}=\frac{1}{4}\sum_k r_k\,d\theta_k.\]
So the squared length of $d\bm{\theta}$ under $g_n$ is
\[\frac{1}{4}\sum_kr_k\,d\theta_k^2-2\biggl( \frac{1}{4}\sum_k r_k\,d\theta_k\!\biggr)^{\!\!2}=
\frac{1}{4}\sum_kr_k\,d\theta_k^2-\frac{1}{8}\Bigl(\sum_k r_k\,d\theta_k\Bigr)^{\!2};\]
and that is $1/8^\text{th}$ of its squared length~(\ref{eq:TangVecBagLength}) under $p_\mathbf{r}$, as we claimed.

\section{Tract boundaries and tract vertices}
\label{sect:VertexPentagons}

Our conformal map of $(E_5, e_5)$ will have the same structure as the repeating icosagon of the Bring sextic in Figure~\ref{fig:Icosagon}; but which equilateral pentagons will lie along the tract boundaries and at the tract vertices on our map?

\subsection{The loops of fixed points}

We first determine the self-inverse, orientation-reversing isometries of $(E_5,e_5)$, since their loops of fixed points are the short-leg and hypot-long-leg loops on~$E_5$.  As happens also on the Bring sextic, there are 26 such isometries.  Just reflecting the pentagons is one such, and we denote that one $\,\flect$; but it has no fixed points.  Indeed, reflecting an $n$\_gon can take it to itself only when that $n$\_gon lies entirely along a single line, and there are no collinear equilateral pentagons.

\begin{figure}
  \begin{center}
	\includegraphics[scale=1]{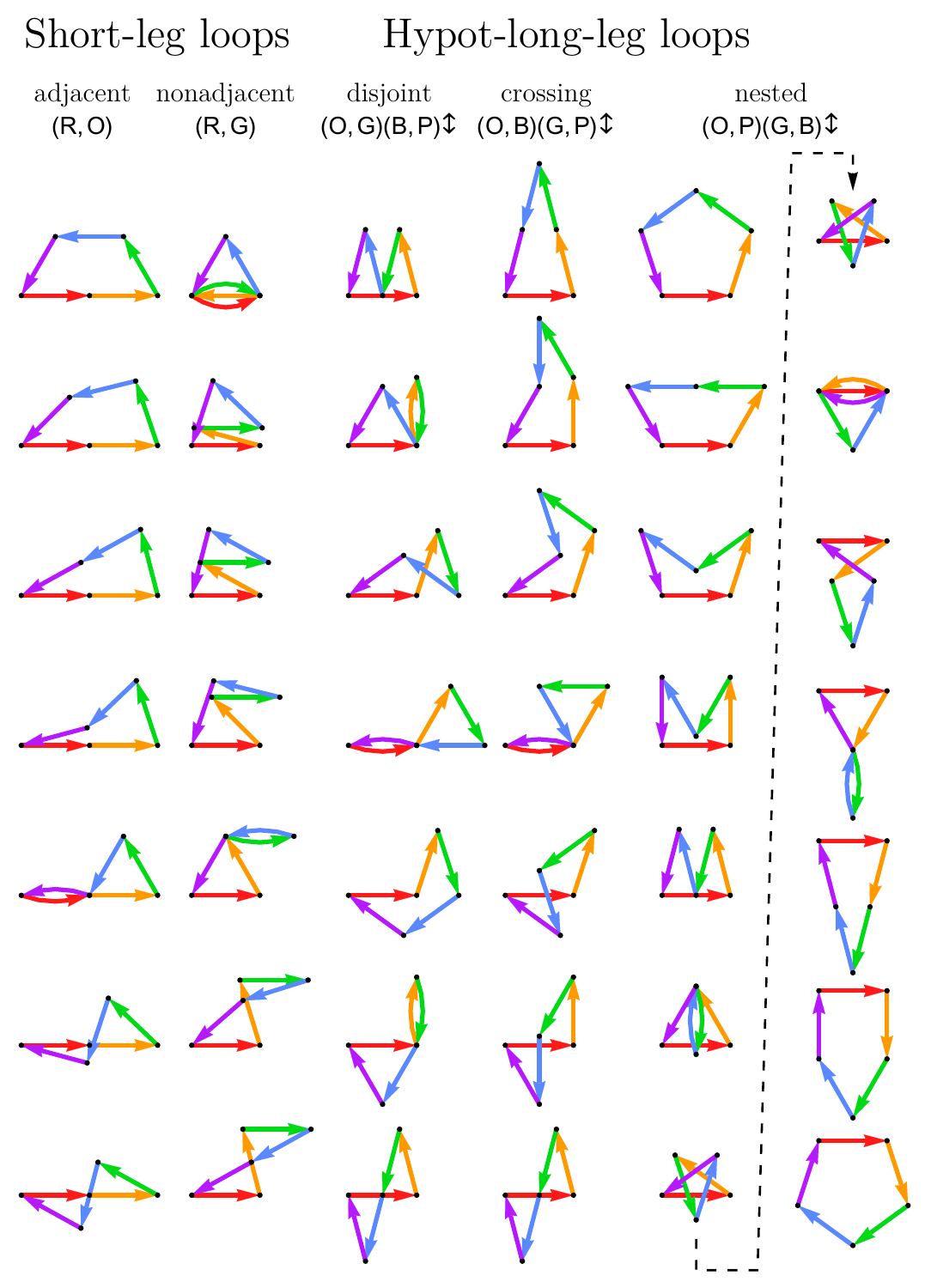}
  \end{center}
	\vspace{-10pt}
	\caption{Sample pentagons from five fixed-point loops on $E_5$, one of each type.  Each of the first four columns depicts a quarter of a loop.  For the next quarter, take the same pictures in the reverse order, but rotate each by $180^\circ$, reverse its arrowheads, and recolor its edges so as to restore the bottom pinwheel, which becomes the joint between the first and second quarters.  That trick doesn\t work for the fifth loop, which doesn\t pass through any pinwheels, so the last two columns together depict half of that fifth loop.  For the second half of any of the five loops, reflect the pentagons of its first half without changing their order.}
	\label{fig:Loops}
\end{figure}

Consider the isometry that just swaps two phases, say the adjacent phases red and orange; we denote that isometry $(\cR,\cO)$, abbreviating the color names.  At least locally on $E_5$, there is a curve of pentagons that remain fixed because their red and orange phases coincide: the $y$\_axis in Figure~\ref{fig:Glimpse}.  Since red is smaller than orange on one side of that curve ($\xi>0$) and larger on the other side ($\xi<0$), the isometry $(\cR,\cO)$ reverses the orientation of $E_5$ by interchanging those two sides.  The pentagons fixed by $(\cR,\cO)$ turn out to form a single loop, a quarter of which is shown in the left column of Figure~\ref{fig:Loops}.  The second column shows a quarter of the loop fixed by swapping the nonadjacent phases red and green: $(\cR,\cG)$.  We call these isometries the \emph{short swaps}, since each fixes a short-leg loop on $E_5$.   There are ten short swaps, five of adjacent pairs and five of nonadjacents.

There is a difference here between our coordinate systems for $E_5$ and for the Bring sextic~$\Bsextic$.  On $\Bsextic$, an antiautomorphism that fixes a short-leg loop swaps some two homogeneous coordinates and also conjugates all five of them.  A short swap on~$E_5$ swaps two phases, but it doesn\t reflect the pentagons, which would negate all five phases (up to an additive constant).  This difference arises from how the coordinates are constrained.  Each point on $\Bsextic$ has complex homogeneous coordinates $(c_1,\ldots,c_5)$, and we can view a pentagon\s edges $(\edge_1,\ldots,\edge_5)$ as being complex homogeneous coordinates for that pentagon.  Both sets of coordinates must sum to zero, which reduces the degrees of freedom from $8$ to $6$.  But the last four degrees of freedom get removed differently:  The $c_k$ also satisfy the two complex equations $\sum_k c_k^2=\sum_k c_k^3=0$, while the five edges $\edge_k$ are constrained to have equal lengths by the four real equations $\edge_1\bar\edge_1=\cdots=\edge_5\bar\edge_5$.  It is those length constraints --- particularly their complex conjugations --- that make swapping a pair of coordinates behave differently on $E_5$ than on $\Bsextic$.

Swapping two disjoint pairs of coordinates preserves the orientation, both on~$\Bsextic$ and on~$E_5$; so there\s no difference here.  To fix a hypot-long-leg loop on $E_5$, we must also reflect the pentagons.  The resulting isometries are the \emph{hypot-long swaps}.  Figure~\ref{fig:Loops} shows a quarter of the loop fixed by $(\cO,\cG)(\cB,\cP)\flect$; a quarter of the loop fixed by $(\cO,\cB)(\cG,\cP)\flect$; and half of the loop fixed by $(\cO,\cP)(\cG,\cB)\flect$.  In the first of those three cases, the swapped pairs of phases are disjoint in the cyclic assembly order of the edges; in the second, they cross each other; in the third, they are nested.  There are $15$ hypot-long swaps, five of each type, and the three types have different geometric implications for the pentagons that they fix:
\begin{describe}
\item[disjoint:] The vertex opposite the unswapped edge lies somewhere on the line extending that edge.
\item[crossing:] If we define an \emph{arm} to consist of two adjacent pentagon edges, the arm entering the vertex opposite the unswapped edge and the arm leaving that vertex differ only by a rigid rotation around that vertex.  
\item[nested:] The pentagon has mirror symmetry across the perpendicular bisector of the unswapped edge (which passes through the vertex opposite that edge).  
\end{describe}

\subsection{The pentagons at the tract vertices}

\begin{figure}
  \begin{center}
	\includegraphics[scale=1.1]{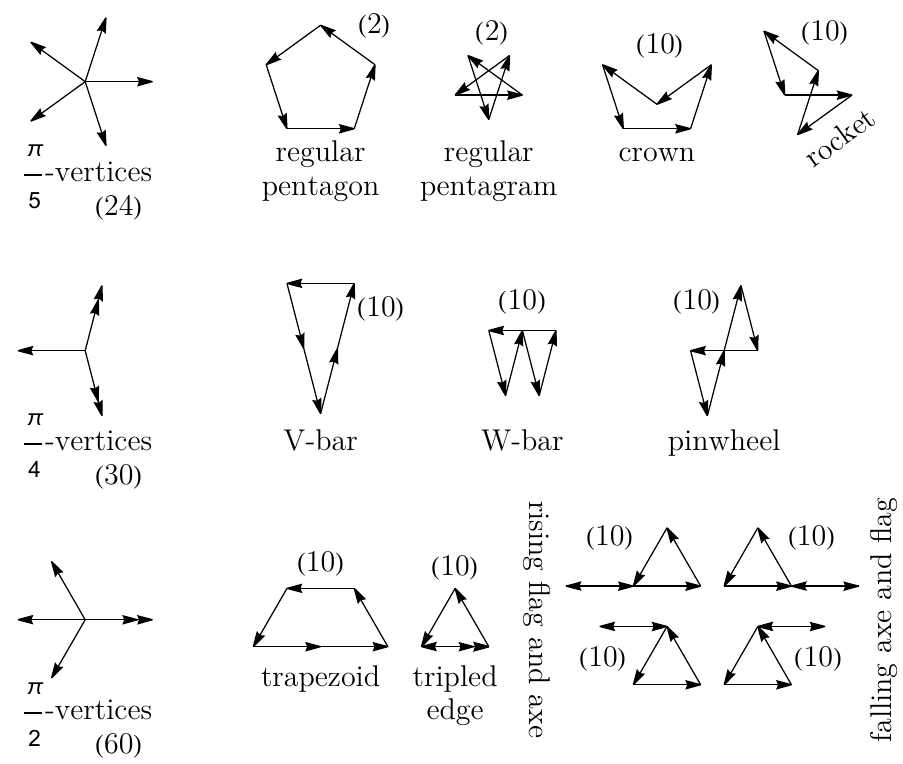}
  \end{center}
	\vspace{-10pt}
	\caption{The pentagons at the $114$ tract vertices of $E_5$ have these $13$ shapes, where the numbers in parentheses say how many instances there are of each shape.  The doubled edge of a \emph{flag} or \emph{axe} is its \emph{handle}, and a flag or axe is \emph{rising} or \emph{falling} according as the edge of its equilateral triangle parallel to that handle points away from or toward the handle.}
	\label{fig:VertexPentagons}
\end{figure}

An equilateral pentagon at a tract vertex lies on more than one of those $25$ loops of fixed points, so its phases are severely constrained.   Figure~\ref{fig:GlimpseAng} has examples of how those constraints can be met, and Figure~\ref{fig:VertexPentagons} catalogs the resulting shapes.
\begin{describe}
\item[$\mathbf{\frac{\pmb{\pi}}{5}}$\_vertices:] The pentagons at the $\frac{\pi}{5}$\_vertices, of which there are $240/10=24$, have their five phases arranged like the complex fifth roots of unity.  Each such vertex lies on five hypot-long-leg loops.
\item[$\mathbf{\frac{\pmb{\pi}}{4}}$\_vertices:] The pentagons at the $\frac{\pi}{4}$\_vertices, of which there are $240/8=30$, have two pairs of coincident phases. Each lies on four loops: two short-leg loops and two hypot-long-leg loops.
\item[$\mathbf{\frac{\pmb{\pi}}{2}}$\_vertices:] The pentagons at the $\frac{\pi}{2}$\_vertices, of which there are $240/4=60$, have their phases arranged like the complex square roots and cube roots of unity superposed. Each lies on one short-leg loop and one hypot-long-leg loop.
\end{describe}

Optionally reflecting each named pentagon in Figure~\ref{fig:VertexPentagons} and then choosing which edge to color red gives ten instances of that named shape --- the exceptions being the regular pentagon and pentagram, whose rotational symmetry makes all color choices equivalent, giving each of them just two instances.  

We distinguish a pinwheel from its reflection by giving the turns, either left or right, \texttt{L} or \texttt{R}, that precede and follow its straight vertex \texttt{S}.  So the pinwheel depicted in Figure~\ref{fig:VertexPentagons} is an \LSR\ pinwheel, while its reflection would be \RSL.  All of the other pentagons depicted in Figure~\ref{fig:VertexPentagons} have positive signed area, so we distinguish them from their reflections by calling them \ccWise; their reflections would be~\cWise.  (Note that the belly of a rocket has more area than its nozzle.)

\section{Sketching our conformal map}
\label{sect:MapSketch}

Suppose that we know which equilateral pentagon lies at each of the three vertices of some tract~$U\mkern -2mu$.  We then know two different pentagons along each boundary of~$U\mkern -2mu$, so we can determine the loop to which each boundary curve belongs.  If $V$ is the tract across one of those curves from $U\mkern -2mu$, we can apply the phase swap of the loop that divides them to the third vertex of~$U\mkern -2mu$, thereby computing the pentagon that lies at the third vertex of~$V\mkern -2mu$.  By repeating that process, we can extend our vertex labeling as far as we like, even to all of $E_5$.

\subsection{One precinct}
\label{subsect:OnePrecinct}

\begin{figure}
  \begin{center}
    \vspace{-5 pt}
	\includegraphics[scale=0.8]{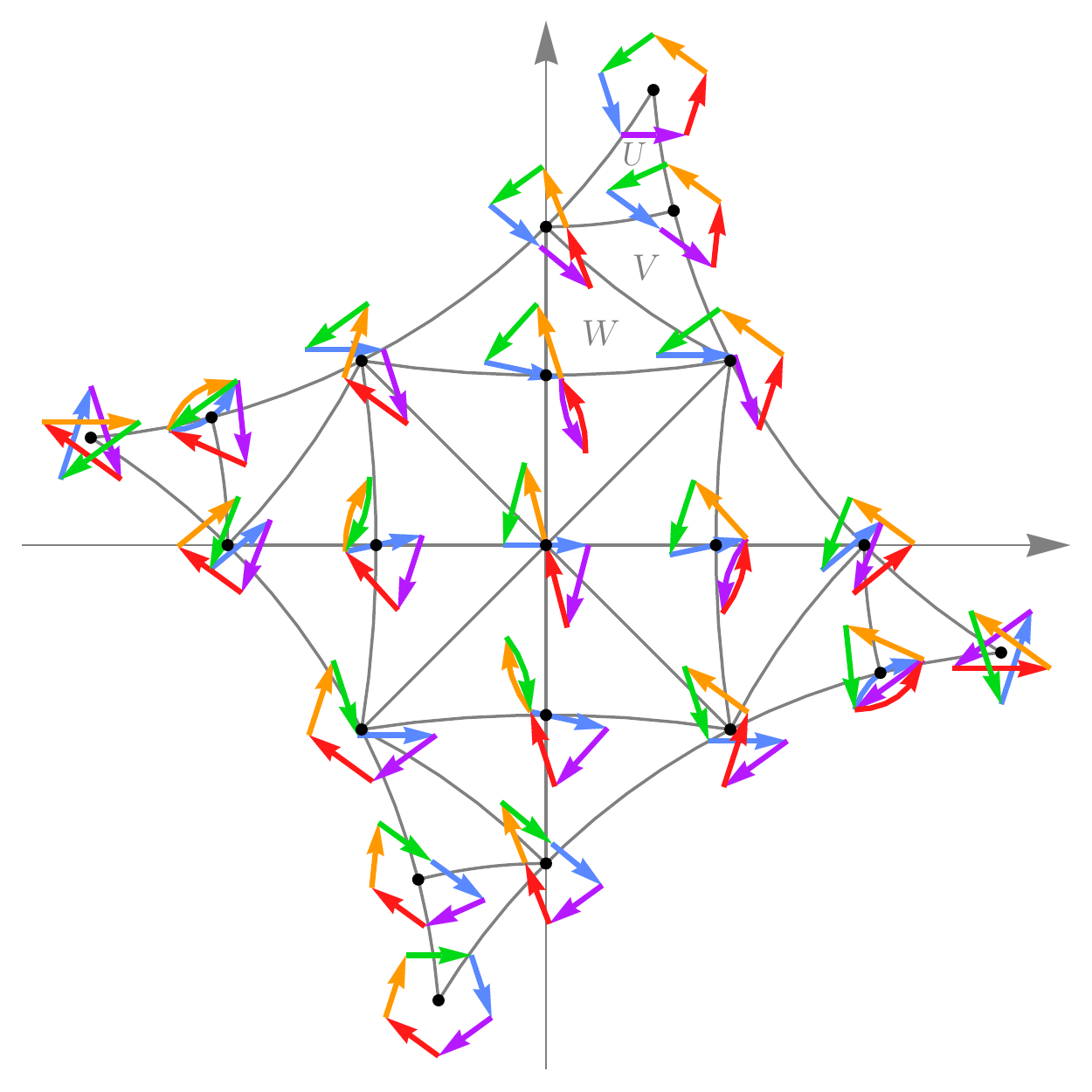}
  \end{center}
	\vspace{-15pt}
	\caption{A precinct with \Rtract\ corners and an \RSL\ pinwheel at its center.}
	\label{fig:Rprecinct}
\end{figure}

In Figure~\ref{fig:Rprecinct}, we have used that process to label the vertices of one precinct.  (Recall, from Section~\ref{sect:BringSextic}, that a \emph{precinct} is a hyperbolic regular quadrilateral with vertex angles of $\frac{\pi}{5}$ that is tiled by $24$ tracts.)  The uppermost tract in that precinct, labeled~$U\mkern -2mu$, is the tract in Figure~\ref{fig:Glimpse} that shares its hypotenuse with $T_\text{out}$, and the rest of the tract vertices have been labeled to be consistent with those of $U\mkern -2mu$.  For example, the short leg of $U$ belongs to the $(\cB,\cP)$ short-leg loop; so the $\frac{\pi}{5}$\_vertex of the southern neighbor $V$ of $U$ is the \ccWise\ crown with an orange base.  

What about the tract $W\mkern -2mu$, the southwest neighbor of $V$?  In the pentagons at both ends of the hypotenuse of~$V\mkern -2mu$, the blue-purple two-edge arm differs from the red-orange arm just by a rigid rotation around the purple-to-red vertex.  So that hypotenuse is part of a hypot-long-leg loop of the crossing type: in particular, the one fixed by the hypot-long swap $(\cR,\cB)(\cO,\cP)\flect$.  Applying that swap to the \ccWise\ trapezoid at the $\frac{\pi}{2}$\_vertex of $V\mkern -2mu$, we find that the pentagon at the $\frac{\pi}{2}$\_vertex of $W$ is the \ccWise\ rising flag with a purple-red handle.  And so forth.

Note that the four hypot-long-leg loops that bound the precinct in Figure~\ref{fig:Rprecinct} are all of the nested type, their unpaired edges having the four colors other than the unpaired blue of the central pinwheel.  We are about to carve up $E_5$ into ten precincts by cutting it along the five nested hypot-long-leg loops.

The precinct in Figure~\ref{fig:Rprecinct} has an \RSL\ pinwheel at its center and \Rtract{s} at all four of its corners.  Cyclically permuting the colors will give us five such precincts.  Reflecting each of those five --- reflecting both the labeling pentagons and the layout of the precinct itself --- will give us five more precincts, each with an \LSR\ pinwheel at its center and \Ltract\ corners.  And those ten precincts turn out to partition~$E_5$.  But wait a minute:  Why do \RSL\ pinwheels go with \Rtract\ corners and \LSR\ pinwheels with \Ltract\ corners, rather than the reverse?

\subsection{The orientation issue}
\label{subsect:Orientation}
 
While the Bring sextic inherits an orientation from its structure as a complex curve, it is less clear how to orient $E_5$ --- or, more generally, any polygon space~$P_\mathbf{r}$ of $n$\_gons.  When $n$ is odd, though, here is one way to choose an orientation.

In the quotient $(S^1)^n/SO(2)$, consider a nonsingular point $\bm{\theta}=(\theta_1,\ldots,\theta_n)$ that belongs to~$P_\mathbf{r}$ because $\sum_k r_ke^{\imagunit\theta_k}=0$.  The tangent space to $(S^1)^n$ at $\bm{\theta}$ is the product of three factors:  the tangent space to $P_\mathbf{r}$ at $\bm{\theta}$ (of dimension $n-3$); the real line~$\reals$, for rigid rotations of the $n$\_gon; and the complex plane $\complexes$, for ways in which the $n$\_gon could fail to close up.  When $n$ is odd, the orientation of that product isn\t affected by the order of its three factors, since only one of them has odd dimension: $\reals$.  The factors $\reals$ and $\complexes$ have standard orientations; and the assembly order of the $n$ edges gives an orientation to~$(S^1)^n$ --- an orientation that, since $n$ is odd, wouldn\t change were that order to be cyclically shifted.  We denote by~$P_\mathbf{r}^+$ the orientation for the first factor that this product implies.

It turns out [A2] that $E_5^+$ is the orientation of $E_5$ in which the precincts with an~\RSL\ pinwheel at their center have \Rtract\ corners, as in Figure~\ref{fig:Rprecinct}.  Indeed, Figure~\ref{fig:Rprecinct} came out that way because this paper adopted the orientation $E_5^+$ from the start.  That\s why, in Figure~\ref{fig:GlimpseAng}, going around from~$\xi$ to $\eta$ passes two other vertices, while going back from~$\eta$ to $\xi$ passes only one, rather than the reverse.

\begin{figure}
  \begin{center}
  \vspace{-10pt}
	\quad\includegraphics[scale=0.9]{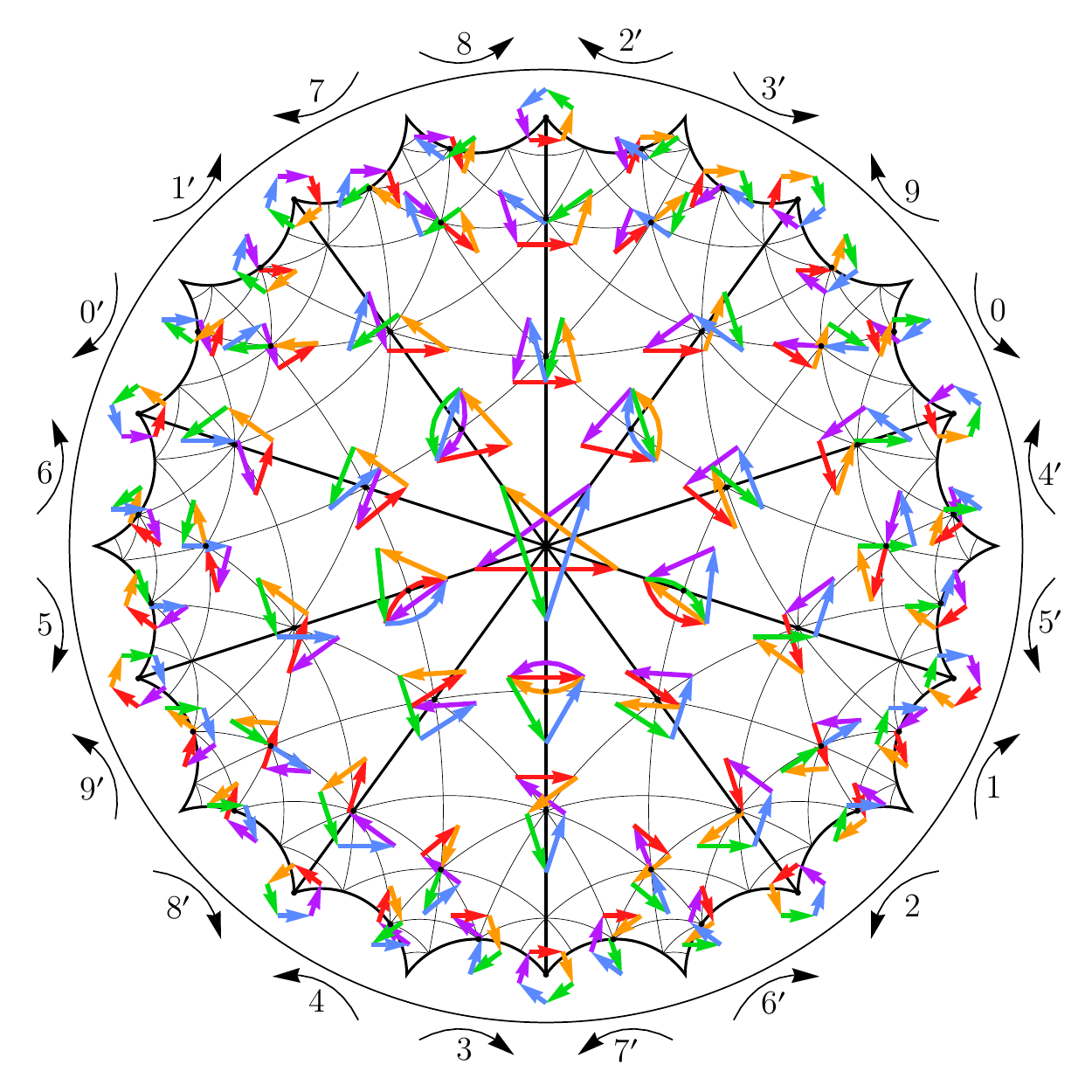}
  \end{center}
	\vspace{-20pt}
	\caption{Our conformal map of the entire polygon space $E_5$, centered at the \ccWise\ regular pentagram and with some tract vertices labeled.  The edges match up as in Figure~\ref{fig:Icosagon}; for example, both~$5$ and~$5'$ contain the \cWise\ crown --- and both~$0$ and~$0'$, the \cWise\ rocket --- on a red base.  (The precinct of Figure~\ref{fig:Rprecinct}, centered at the \RSL\ pinwheel whose unpaired edge is blue, is here at 9 o'clock.)}
	\label{fig:LabeledIcos}
\end{figure}

For an easy test of which orientation $E_5$ has been given in some diagram, examine any tract that has a regular pentagon at its $\frac{\pi}{5}$\_vertex, such as the top or bottom tract in Figure~\ref{fig:Rprecinct}.  Rotate the labeling pentagons along that tract\s short leg so that their double-length sides are parallel.  When the orientation is~$E_5^+$, the edge opposite that double-length side (that edge being orange at the top of Figure~\ref{fig:Rprecinct}, red at the bottom) rotates \cWise\ as we revolve \cWise\ around the regular pentagon at the $\frac{\pi}{5}$\_vertex, and ditto for \ccWise\ --- like the earth, which rotates around its axis in the same direction that it revolves around the sun.\footnote{But it works the other way for a tract of $E_5^+$ with a regular pentagram at its $\frac{\pi}{5}$\_vertex.  In each pentagon along such a short leg, the two edges with the same phase are nonadjacent.  Rotating those pentagons to keep that shared phase constant, the edge between that same-phase pair (which is green at the left in Figure~\ref{fig:Rprecinct}, purple at the right) rotates \cWise\ as we revolve \ccWise\ around the regular pentagram at the $\frac{\pi}{5}$\_vertex, and vice versa.}

\subsection{Assembling the icosagon}

By reflecting and cyclically recoloring Figure~\ref{fig:Rprecinct}, we produce ten precincts that partition~$E_5$.  We can cluster those ten around any one of their four shared corners to form an icosagon as in Figure~\ref{fig:Icosagon}.   Figure~\ref{fig:LabeledIcos} clusters them around the \ccWise\ regular pentagram; so the \cWise\ regular pentagram appears at the ten $\smash{\frac{\pi}{5}}$\_vertices (not labeled, due to the cramped quarters), while the \cWise\ and \ccWise\ regular pentagons alternate at the $\frac{2\pi}{5}$\_vertices.  If we clustered the ten instead around the \ccWise\ regular pentagon, the result would extend Figure~\ref{fig:Glimpse} to cover all of $E_5$.

A warning about how the labeling pentagons are rotated:  When building a polygon space, we identify pentagons that differ by rotation; but we must choose some rotation for each pentagon that labels a map.  When drawing~$E_4$ in Figure~\ref{fig:Tetravecs} and when tracking the loops of $E_5$ in Figure~\ref{fig:Loops}, we kept the red phase consistently at zero.  In most of our drawings of $E_5$, however, to avoid singling out one edge, we keep the sum of the five phases at a fixed value.  But rotating a pentagon by $2\pi/5$ doesn\t change that sum modulo $2\pi$, so this rule leaves five choices for each pentagon; and a trip around $E_5$ can bring us back to the same pentagon with a different rotation.  For example, the five \ccWise\ regular pentagons at the heads of the even-numbered arrows in Figure~\ref{fig:LabeledIcos} are all differently rotated.

\section{The dodecadodecahedral model}
\label{sect:Dodeca}  

The global connectivity of $E_5$ is captured, in principle, by the subtle way in which the icosagon of Figure~\ref{fig:LabeledIcos} repeats as we travel around the hyperbolic plane.  But modeling $E_5$ as a surface in $3$\_space makes that connectivity much clearer.

\subsection{Modeling isometries as Euclidean symmetries}

It\s easy to embed a quadruple torus smoothly in $3$\_space; but we want as many of the isometries of $(E_5,e_5)$ as possible to become Euclidean symmetries of our model.  Permuting the phases through the $5$\_cycle $(\theta_1,\ldots,\theta_5)$ is an orientation-preserving isometry that fixes four points:  the two regular pentagons and the two regular pentagrams.  It rotates the neighborhoods on~$E_5$ of the pentagons by $\pm 4\pi/5$, as shown by cycling the colors in Figure~\ref{fig:Glimpse}.  But it rotates the neighborhoods of the pentagrams by $\pm 2\pi/5$, as shown in Figure~\ref{fig:LabeledIcos}.  Handling that discrepancy in $4$\_space would be no problem.  Indeed, cyclically permuting the dimensions of $(S^1)^5$ effects a double rotation of the $4$\_dimensional quotient $(S^1)^5/SO(2)$, a rotation that carries the subset $E_5$ to itself while rotating it by $\pm 2\pi/5$ around one axis and by $\pm 4\pi/5$ around the other.  But dealing with that discrepancy in $3$\_space is a challenge, since rotations there have only one axis.  Note that doubling $\pm 4\pi/5$ gives $\pm 8\pi/5$, which is congruent modulo $2\pi$ to $\mp 2\pi/5$; thus, when we ignore signs, doubling these angles is the same as halving them.

\begin{figure}[!tb]
	\begin{center}
	\includegraphics[scale=0.45]{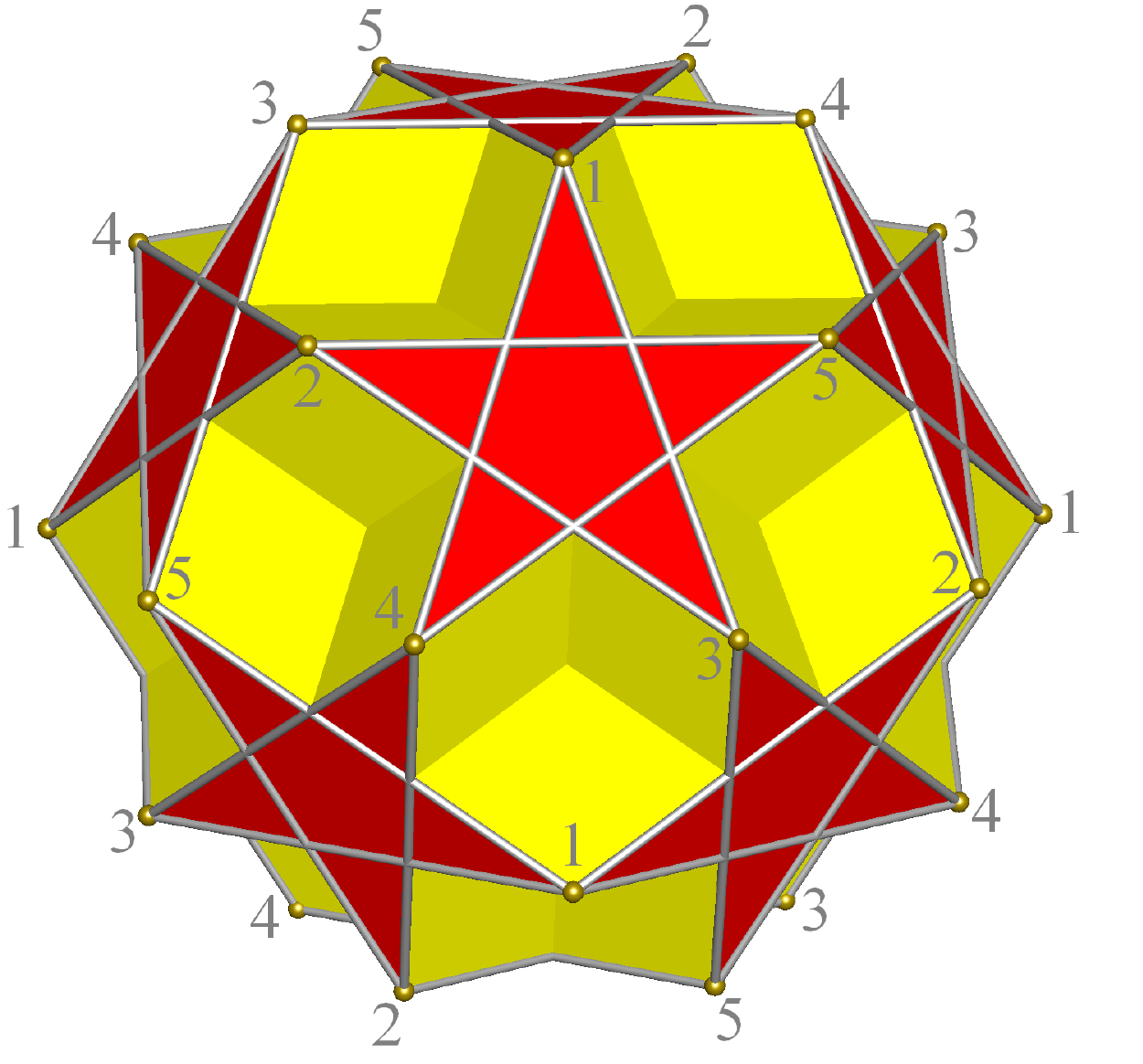}\quad
	\includegraphics[scale=0.44]{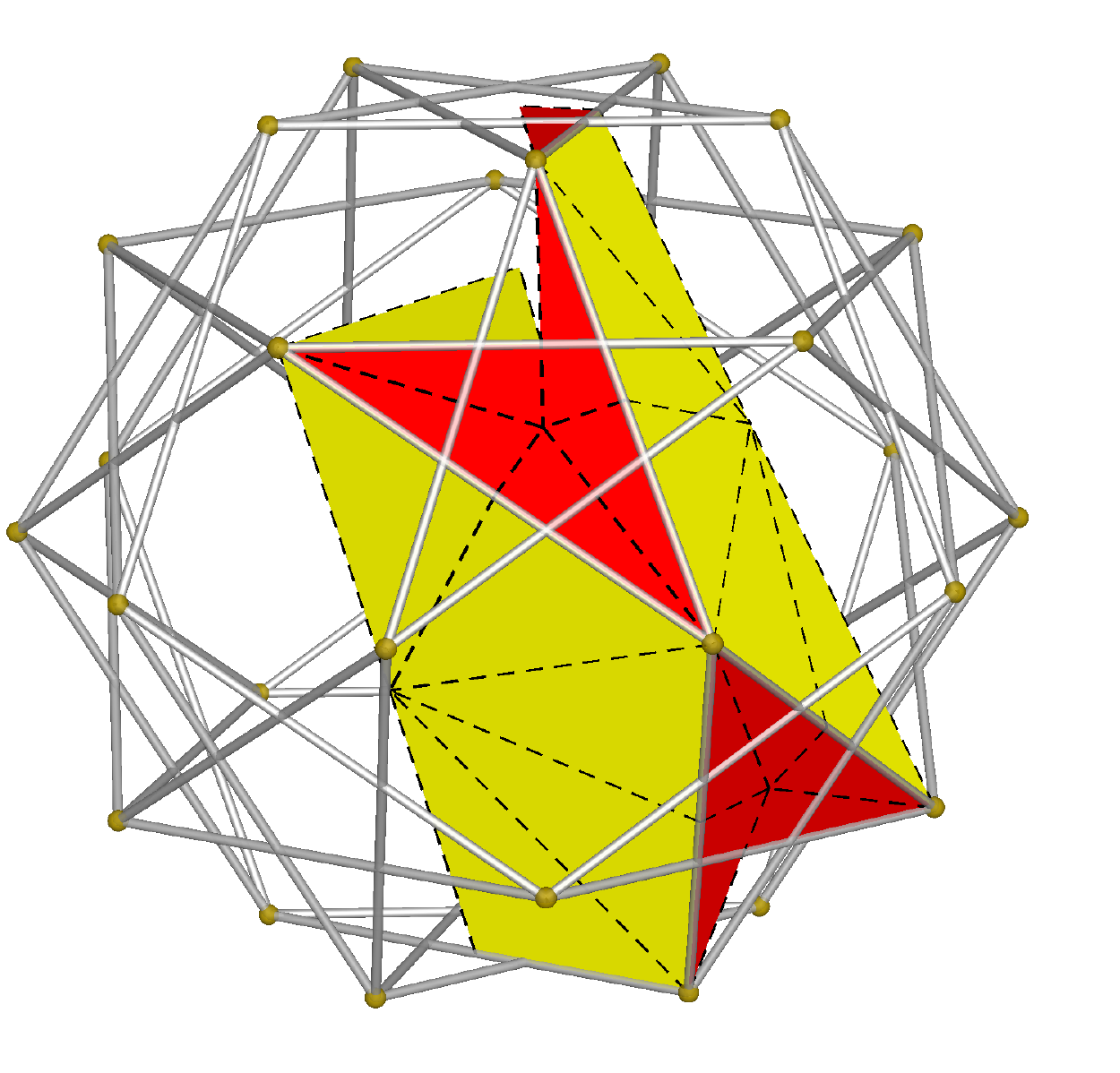}
	\end{center}
	\vspace{-10pt}
	\caption{The left shows a dodecadodecahedron, drawn in the graphic style of Stella~\cite{Stella}.  The points labeled $k$, for $1\le k\le 5$, are vertices of the $k^\text{th}$ inscribed regular octahedron.  The right shows the $24$ tracts of one of the ten precincts into which $E_5$ is partitioned, like the one in Figure~\ref{fig:Rprecinct}, assuming that the pinwheel pentagons are at the ten equatorial vertices.  (Double wrapping makes the first and fifth tracts in a half of a red decatract appear to bump into each other.)}
	\label{fig:Dodeca}
\end{figure}

\begin{figure}
    \vspace{5 pt}
	\begin{center}
	\includegraphics[scale=0.5]{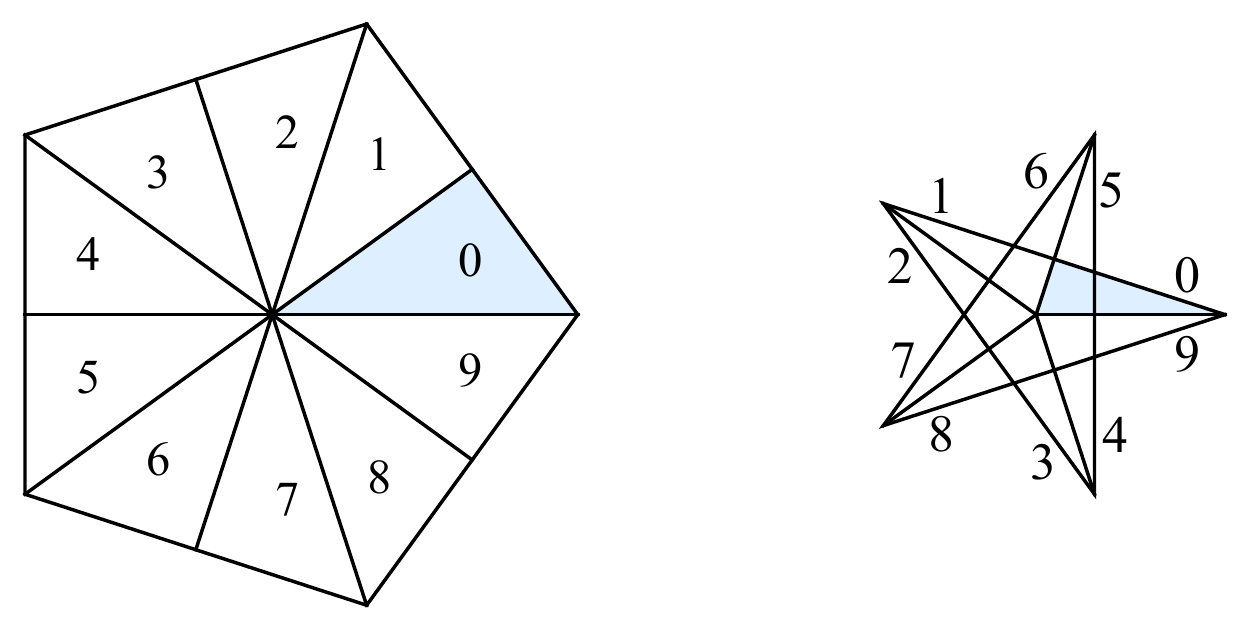}
	\end{center}
	\vspace{-15pt}
	\caption{How Weber models a decatract using either a pentagon or, via double wrapping, a pentagram.  Both of the blue-tinted tracts are \Rtract{s}, even though the long (radial) leg of a tract on the right is shorter than its short leg.}
	\label{fig:DoubleWrapping}
\end{figure}

Weber~\cite{Weber} suggests modeling the Bring sextic $\Bsextic$ using a \emph{dodecadodecahedron}, the self-intersecting polyhedron shown in Figure~\ref{fig:Dodeca}.  It has, as faces, $12$ regular pentagons and $12$ regular pentagrams.  Given the $24$ decatracts of $\Bsextic$ (or $E_5$), Weber models one on each pentagon in the obvious way and one on each pentagram by wrapping it twice around its center in the way that the squaring function $z\mapsto z^2$ doubly wraps the neighborhood of the origin in~$\complexes$, as shown in Figure~\ref{fig:DoubleWrapping}.  This handles that discrepancy:  Rotating the model by any multiple of $2\pi/5$, the decatracts on the orthogonal pentagrams rotate, ignoring signs, both twice as far and half as far as those on the orthogonal pentagons.  Introducing those $12$ double points also reduces the Euler characteristic of the model from the $+6$ of three nested spheres to the $-6$ that is required for $\Bsextic$ (and $E_5$).

\subsection{The checkerboard structures}

The dodecadodecahedron\s $24$ faces alternate in a checkerboard fashion between pentagons and pentagrams.  The decatracts of $E_5$ correspond to the $24$ cyclic orders of the five edges by phase, and that order flips between even and odd each time that we enter a new decatract by crossing a short leg; so the decatracts of~$E_5$ also have a checkerboard structure.\footnote{The decatracts of $\Bsextic$ correspond to the $24$ cyclic orders of the five homogeneous coordinates by phase, and that order also flips between even and odd when entering a new decatract; so~$\Bsextic$ has an analogous checkerboard structure.  But assembling five edge vectors in different orders gives dramatically different pentagons, while indexing five homogeneous-coordinate values in different orders has little effect; so the checkerboard structure on $E_5$ matters more than the one on $\Bsextic$.}  Modeling $E_5$ as a dodecadodecahedron requires aligning those two checkerboards.  If we put the decatracts centered at regular pentagons or rockets onto pentagonal faces, then those centered at regular pentagrams or crowns end up on pentagrammatic faces, which seems attractive.

Perhaps the prettiest choice uses the upper horizontal pentagonal face to model the decatract of convex-and-\ccWise\ pentagons shown in Figure~\ref{fig:Glimpse}, with the \ccWise\ regular pentagon at its center.  The \ccWise\ regular pentagram then centers the north-polar pentagrammatic face.  The five \ccWise\ rockets center the north-temperate pentagonal faces, while the five \ccWise\ crowns center the north-temperate pentagrammatic faces.  (The $20$ slanted faces form five parallel quadruples; and the four faces in each quadruple have, at their centers, the \cWise\ and \ccWise\ crowns and rockets whose bases share a common color.)  The five north-temperate vertices are the \ccWise\ V\_bars, the five north-polar vertices are the \ccWise\ W\_bars, and the ten equatorial vertices are the pinwheels, alternating between \RSL\footnote{Assuming the orientation $E_5^+$, the front equatorial vertex numbered 3 in Figure~\ref{fig:Dodeca} is an \RSL\ pinwheel, since it lies at the center of a precinct, shown on the right, that has \Rtract\ corners.} and~\LSR.  The \ccWise\ $\frac{\pi}{2}$\_vertices --- five each of trapezoids, tripled edges, rising flags, falling flags, rising axes, and falling axes --- bisect the $30$ edges that lie in the northern hemisphere.

The southern hemisphere is entirely analogous, but all~\cWise.  The pentagons with zero signed area form an ``equator'' that loops around the tropics three times.  It passes once through each pinwheel pentagon, but not in longitude order; and it travels across each of the ten  rocket-centered, slanted pentagonal faces in order to connect that face\s two pinwheel vertices.\footnote{The locus of zero-signed-area pentagons is a long way from being a geodesic under either the bag-of-edges metric $e_5$ or its hyperbolic rescaling.  As that locus passes through one of the pinwheel vertices of a rocket-centered decatract, the angle that it makes with the decatract edge is $\arctan(3/4)\approx 36.87^\circ$.   In contrast, the bag-of-edges geodesic joining the two pinwheel vertices of that decatract makes an angle with the decatract edge of only about $30.89^\circ$, while that angle for the nearby hyperbolic geodesic is $\arccot(\phi)\approx 31.72^\circ$, where $\phi=(\sqrt{5}+1)/2$ is the golden ratio.  (Rescaling thus twists this geodesic by less than $1^\circ$, as is true also for the similar geodesic in the upcoming Section~\ref{subsect:Ditri}.)  Extending either geodesic would give a geodesic loop whose removal from~$E_5$ would disconnect it, the hyperbolic one having length $10\arccosh(\phi^2)\approx 16.169$.}

\subsection{Even versus odd permutations of the phases}

A dodecadodecahedron has five inscribed regular octahedra, and its $60$ rotations permute them evenly.  Its $30$ vertices, six in each octahedron, model the $\frac{\pi}{4}$\_vertices of $E_5$.  At each $\frac{\pi}{4}$\_vertex, two pairs of phases coincide; and the octahedron to which such a vertex belongs is determined by its unpaired phase.  The rotations of the model thus evenly permute the five phases also.

The $120$ icosahedral symmetries of the dodecadodecahedron consist of the $60$ rotations, which evenly permute the octahedra (on $E_5$, the phases), and each of those combined with central inversion (on $E_5$, with reflection\footnote{Figuring out the complicated way in which reflection transforms the repeating icosagon in Figure~\ref{fig:LabeledIcos} makes me appreciate the dodecadodecahedral model all the more.}).  Together, these $120$ symmetries form the group $A_5\times C_2$.  The group formed by the $240$ isometries of $(E_5,e_5)$ is $S_5\times C_2$.  Since neither~$S_5\times C_2$ nor its subgroup $S_5$ occurs as a point group in three dimensions, its subgroup $A_5\times C_2$ is the only way to realize even half of the isometries of $E_5$ as Euclidean symmetries of a model in $3$\_space.  

An isometry of $E_5$ that effects an odd permutation of the phases interchanges the pentagonal faces with the pentagrammatic ones.  As a taste of the various resulting geometries, let\s consider a short swap.  The short-leg loop where the two swapped phases coincide is a great hexagon of the dodecadodecahedron (closer to horizontal for an adjacent short swap, closer to vertical for a nonadjacent one).  That hexagon remains pointwise fixed while the pentagonal and pentagrammatic faces on the two sides of each of its six edges, say $F$ and $F'$, get interchanged.  Each face $G$ that isn\t adjacent to the fixed hexagon is parallel to two faces that are adjacent to it, say $F_1$ and $F_2$, where $F_1$ and $F_2$ share the wrapping number opposite to that of $G$.  The face $G$ interchanges with the unique face $G'$ that is parallel to $F'_1$ and $F'_2$ but distinct from them and not adjacent to $G$.

\section{Using vertex angles as coordinates}
\label{sect:AngleCoords}

We have come as far as the easy arguments can take us.  To compute the full pentagon-to-point correspondence of our conformal map, we have to work harder.  That work starts by studying the vertex angles of an equilateral pentagon.

\subsection{Relations among the vertex angles}

How are the five vertex angles of an equilateral pentagon related?  The vertex angles are the differences of adjacent phases; we denote them $\tau_k=\theta_{k+1}-\theta_k$, treating the subscripts modulo $5$.  For a convex-and-\ccWise\ pentagon, the $\tau_k$ are the external angles at the vertices, so they are nonnegative and sum to $2\pi$.\footnote{If $\alpha$, $\beta$, $\gamma$, and $\delta$ are the external angles at four successive vertices of a convex-and-\ccWise\ equilateral pentagon, an exercise in Euclidean geometry shows that $\alpha\ge\delta$ just when $\beta\le\gamma$.  It follows that, of the five vertex angles of such a pentagon, the largest always has the two smallest as its two neighbors and the smallest has the two largest as its neighbors, leaving the medial one to fill the remaining slot.  We noted the adjacency of the smallest and largest in Figure~\ref{fig:GlimpseAng}.}  In general, though, the $\tau_k$ are defined only modulo $2\pi$.

To study how the $\tau_k$ are related, we rotate the pentagon so that $\theta_3=0$.  We then have $\theta_4=\tau_3$, $\theta_5=\tau_3+\tau_4$, $\theta_2=-\tau_2$, and $\theta_1=-\tau_2-\tau_1$.  Since $\sum_k\cos\theta_k=\sum_k\sin\theta_k=0$, it follows that
\begin{align}
\notag
\!1+\cos\tau_3+\cos\tau_3\cos\tau_4-\sin\tau_3\sin\tau_4+
\cos\tau_2+\cos\tau_2\cos\tau_1-\sin\tau_2\sin\tau_1&=0\quad\\[-8pt]
\label{eq:trigsums}
&\\[-6pt]
\notag
\sin\tau_3+\sin\tau_3\cos\tau_4+\cos\tau_3\sin\tau_4-
\sin\tau_2-\sin\tau_2\cos\tau_1-\cos\tau_2\sin\tau_1&=0.\quad
\end{align}

We retreat from trigonometry to algebra via the Weierstrass transformation $t_k=\tan(\tau_k/2)$, substituting $\cos\tau_k=(1-t_k^2)/(1+t_k^2)$ or $\sin\tau_k=2t_k/(1+t_k^2)$ for each cosine or sine in~(\ref{eq:trigsums}).  Multiplying through by the common denominator of $(1+t_1^2)\cdots(1+t_4^2)$, which is positive, we get $24$\_term polynomials $Z_\text{cos}(t_1,t_2,t_3,t_4)$ and $Z_\text{sin}(t_1,t_2,t_3,t_4)$ that must be zero [A3].

The ideal generated by $Z_\text{cos}$ and $Z_\text{sin}$ contains polynomials that relate any three of the four variables --- trivariate polynomials that we can compute via Gr\"obner bases.  We get two such trivariates by eliminating first $t_4$ and then $t_2$ [A3]:
\begin{align}
\label{eq:ThreeAdj}
\begin{split}
15+3t_1^2-t_2^2+3t_3^2-16t_1t_2-8t_1t_3-16t_2t_3\hskip 1.3in\\[-2 pt]
\null +3t_1^2t_2^2-t_1^2t_3^2+3t_2^2t_3^2+8t_1t_2^2t_3-t_1^2t_2^2t_3^2&=0
\end{split}\\[4 pt]
\label{eq:ThreeNonAdj}
5+9t_1^2-3t_3^2-8t_3t_4-3t_4^2+t_1^2t_3^2-8t_1^2t_3t_4+t_1^2t_4^2-3t_3^2t_4^2+
t_1^2t_3^2t_4^2&=0.
\end{align}
With the subscripts cyclically shifted, (\ref{eq:ThreeAdj}) relates any three adjacent vertices of an equilateral pentagon, while (\ref{eq:ThreeNonAdj}) handles the nonadjacent triples, relating a vertex angle (such as~$\tau_1$) to the two angles to which it is not adjacent ($\tau_3$ and~$\tau_4$). 

Reflecting an equilateral pentagon negates the $\tau_k$, which negates the~$t_k$; this preserves (\ref{eq:ThreeAdj}) and (\ref{eq:ThreeNonAdj}), since all of their terms have even total degree.  Reindexing the phases in reverse order also preserves (\ref{eq:ThreeAdj}) and (\ref{eq:ThreeNonAdj}), since~$t_1$ and~$t_3$ appear symmetrically in (\ref{eq:ThreeAdj}), while~$t_3$ and~$t_4$ appear symmetrically in (\ref{eq:ThreeNonAdj}). (By the way, we use only (\ref{eq:ThreeNonAdj}) below, having computed (\ref{eq:ThreeAdj}) just for completeness.)

\subsection{Our naive coordinate system}

We choose the vertex angles $\tau_1$ and $\tau_4$ as naive coordinates on $E_5$, giving them the new names $\xi$ and $\eta$.  This is consistent with Figure~\ref{fig:GlimpseAng}, where the five edges $(\edge_1,\ldots,\edge_5)$, with phases $(\theta_1,\ldots,\theta_5)$, are colored red, orange, green, blue, and purple, so that $\xi=\tau_1=\theta_2-\theta_1$ and $\eta=\tau_4=\theta_5-\theta_4$ are the angles at the red-to-orange and blue-to-purple vertices.  We also introduce $x=\tan(\xi/2)$ and $y=\tan(\eta/2)$ as new names for the Weierstrass transforms $t_1$ and $t_4$.

More abstractly, we are taking the phase-difference model $E_5^\Delta$ of $E_5$ in the all-pairs torus $(S^1)^{10}$ and projecting it orthographically down onto the coordinate plane --- actually, the coordinate $2$\_torus --- whose coordinates are $(\tau_1,\tau_4)=(\xi,\eta)$.

Given values for $\xi$ and $\eta$, what pentagons are possible?  It suffices to determine the vertex angles $\tau_2$ and $\tau_3$, since four vertex angles uniquely determine the fifth.  We can solve (\ref{eq:ThreeNonAdj}) to get two possible values for $t_3$.  Swapping~$x$ with~$y$ replaces $(x,t_2,t_3,y)$ with $(y,t_3,t_2,x)$ and hence gives us the two possible values for~$t_2$.  Plugging the four combinations of $t_2$ and $t_3$ into $Z_\text{cos}$ and $Z_\text{sin}$ reveals [A4] that there are just two possible pentagons:
\begin{equation}
\label{eq:stformulas}
t_2=\frac{4x(1+y^2)\pm\sqrt{Q}}{(1+x^2)(y^2-3)}\qquad\text{and}\qquad
t_3=\frac{4y(1+x^2)\pm\sqrt{Q}}{(1+y^2)(x^2-3)},
\end{equation}
with the same sign in the two numerators, where the radicand $Q$ is
\begin{equation}
\label{eq:subexpr}
Q=15+22x^2+22y^2-9x^4+60x^2y^2-9y^4+6x^4y^2+6x^2y^4-x^4y^4.
\end{equation}

\begin{figure}
	\begin{center}
	\includegraphics[scale=0.4]{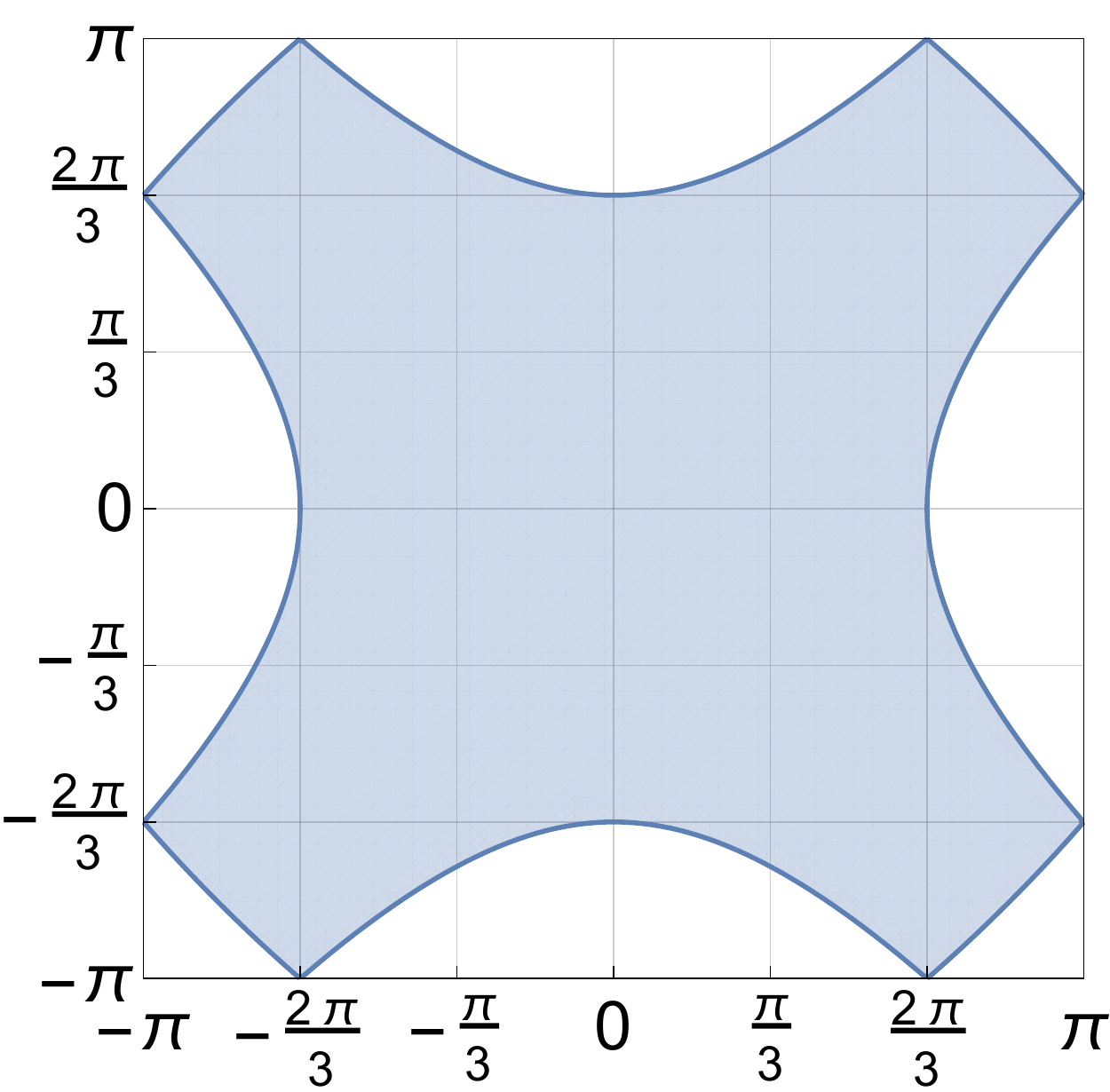} %from TractBoundaries.nb
	\end{center}
	\vspace{-10pt}
	\caption{The region $R$ of the $(\xi,\eta)$ torus where the radicand $Q$ is nonnegative.}
	\label{fig:RegionR}
\end{figure}

Figure~\ref{fig:RegionR} shows the region of our coordinate $2$\_torus where $Q\ge 0$.  Each point $(\xi,\eta)$ where $Q>0$ models two pentagons, the convex-and-\ccWise\ pentagons of Figure~1 arising from the negative signs in~(\ref{eq:stformulas}).  Where $Q=0$, there is only one, and where $Q<0$, none.  So we are modeling $E_5$ as two copies of~$R$: two sheets sewn together along the boundary curve $Q=0$.  If~$R$ were a disk, $E_5$ would be topologically a $2$\_sphere.  But $R$ wraps around the torus both ways, touching itself at the points where $(\xi,\eta)$ is $(\pm 2\pi/3,\pi)$ or $(\pi,\pm 2\pi/3)$.  What\s going on there?

\subsection{Ditri-pentagon loops}
\label{subsect:Ditri}

A pentagon that decomposes must somehow intermix the edges of a digon and a triangle.  In the equilateral case, let\s call such a thing a \emph{ditri-pentagon}.

In Figure~\ref{fig:RegionR}, consider pentagons in which $\eta=\pi$.  The edges $\edge_4$ and $\edge_5$ then form a digon, which forces $\edge_1$, $\edge_2$, and $\edge_3$ to form an equilateral triangle, either \cWise\ or \ccWise.  Since the angle between the digon and the triangle is free to vary, there are two entire loops of such ditri-pentagons.  Our coordinate system collapses each of those loops down into a single point: $(\xi,\eta)=(\pm 2\pi/3,\pi)$.  It also collapses two loops with $\xi=\pi$ into the points $(\xi,\eta)=(\pi,\pm 2\pi/3)$.  Each such collapse pinches off one of the four handles by which $E_5$ differs topologically from a $2$\_sphere.

A pentagon has ten pairs of edges, so there are $20$ ditri-pentagon loops.  As we saw in Section~\ref{sect:PhaseDifference}, each such loop maps to a straight line  in the all-pairs torus $(S^1)^{10}$, a line that lies entirely in the phase-difference model $E_5^\Delta$ of $E_5$.  Our $(\xi,\eta)$ coordinate system collapses four of those lines down into the wrap-around points.  It projects the other $16$ onto the straight contours in the upcoming Figure~\ref{fig:HomeDistrictK}, the ones labeled $-50/21$, four with each of the four slopes $0$, $\infty$, $+1$, and $-1$.

Since the ditri-pentagon loops are straight lines in the phase-difference model of $E_5$, they are geodesics of the Riemannian manifold $(E_5,e_5)$.  The metric $e_5$ must be rescaled a bit, though, to give $E_5$ constant Gaussian curvature, so the ditri-pentagon loops are not hyperbolic straight lines on our conformal map.  Instead, each of the $20$ ditri-pentagon loops snakes gently back and forth across one of the $20$ altitude loops, which is hyperbolically straight.  The two loops in such a pair make an angle of about~$0.95^\circ$ where they cross at each $\frac{\pi}{2}$\_vertex, as shown in Figure~\ref{fig:ScaleFactor}.\footnote{Where the loops cross at a $\frac{\pi}{2}$\_vertex, the angle between the short leg and the ditri-pentagon loop is $\mydepthstrut{4pt}\smash{\arctan\bigl(\sqrt{5/7}\bigr)}\approx 40.20^\circ$, but between short leg and altitude is
$\mystrut{6.5pt}\smash{\arctan\bigl(\sqrt{3-\mystrut{6.5pt}\smash{\sqrt{%
\mystrut{4pt}\smash{5}}}}\,\bigr)}\approx 41.15^\circ$.}  The rescaling stretches the neighborhoods of the $\smash{\frac{\pi}{4}}$\_vertices the most, and the altitude loops respond to that by staying a skosh farther away from the $\smash{\frac{\pi}{4}}$\_vertices, in order to travel through territory that has been stretched less.

\subsection{The tract boundaries}

\begin{figure}
  \vspace*{-10pt}
  \begin{center}\qquad\qquad\qquad\qquad\qquad\quad
    \includegraphics[scale=0.4]{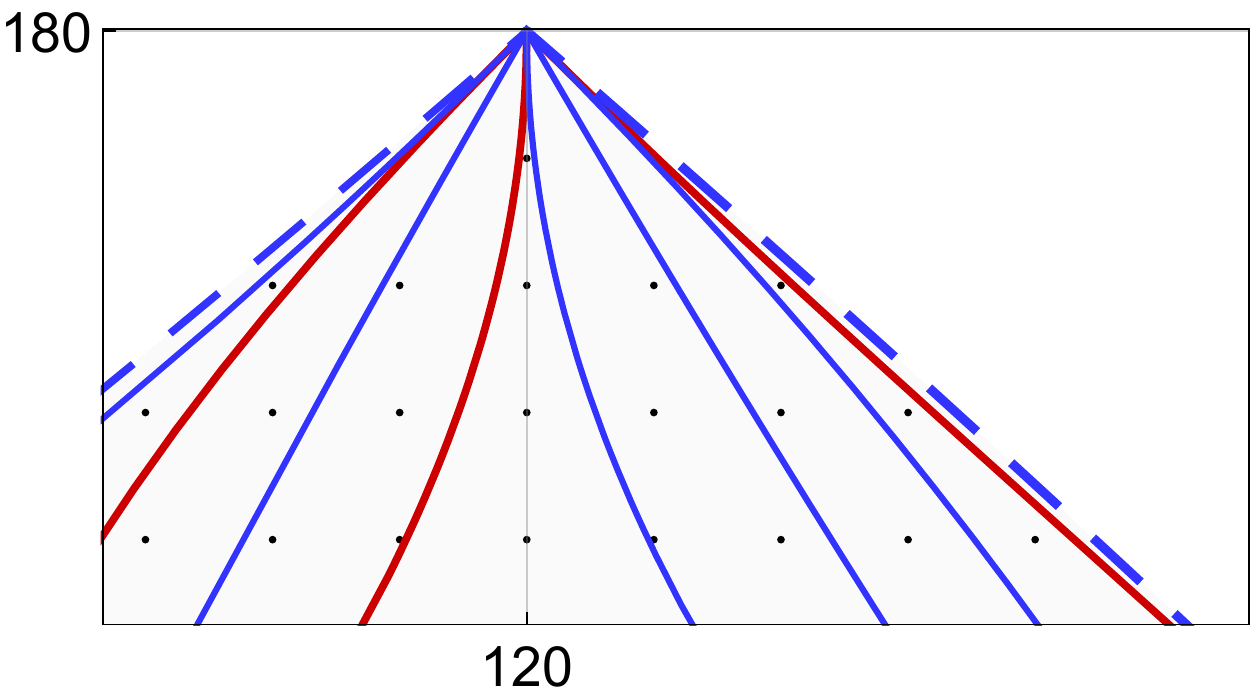} % from TractBoundaries.nb
    \end{center}
    \vspace*{-20pt}
    \begin{center}
	\includegraphics[scale=0.62]{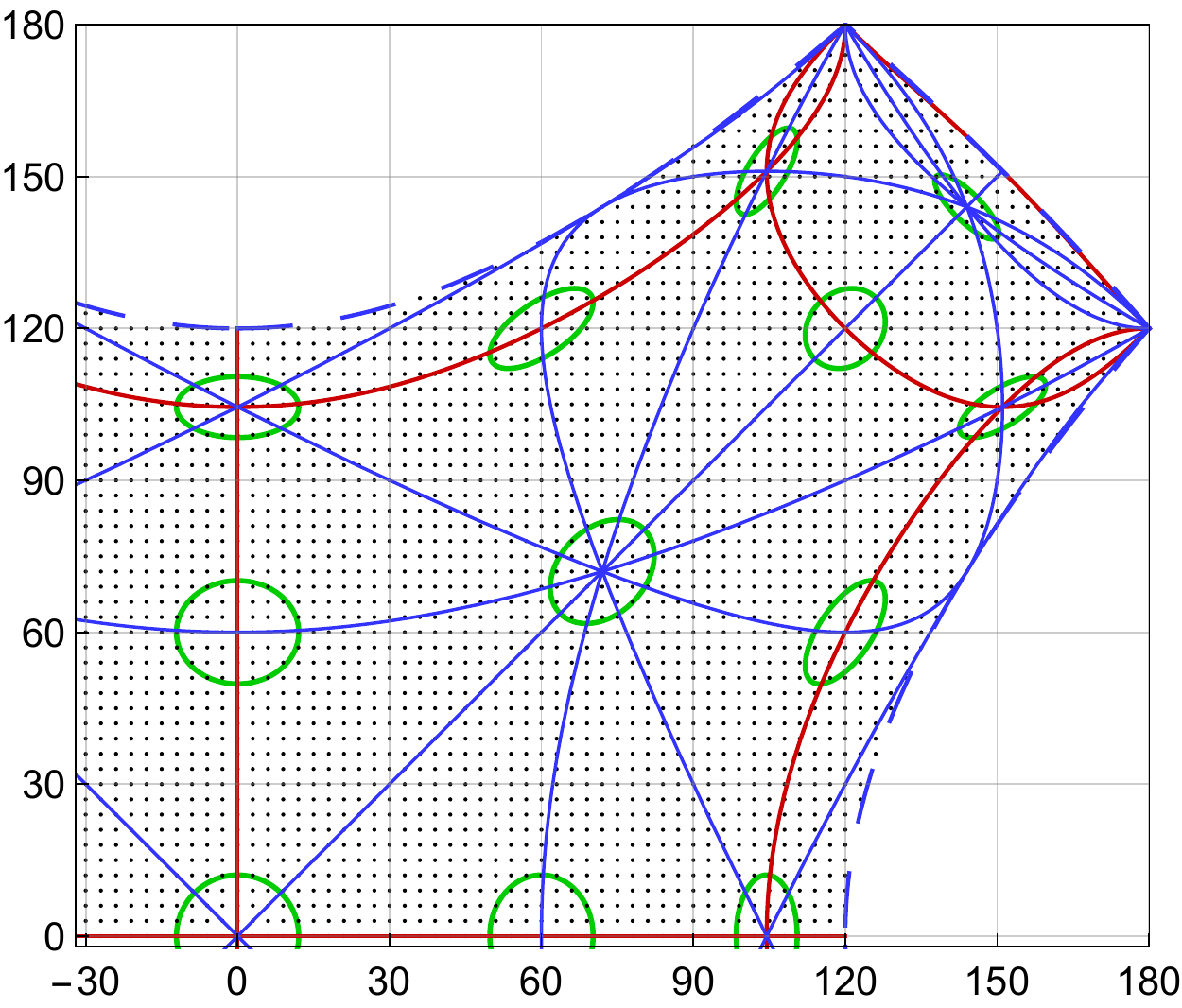}
  \end{center}
  \begin{center}
  \includegraphics[scale=0.58]{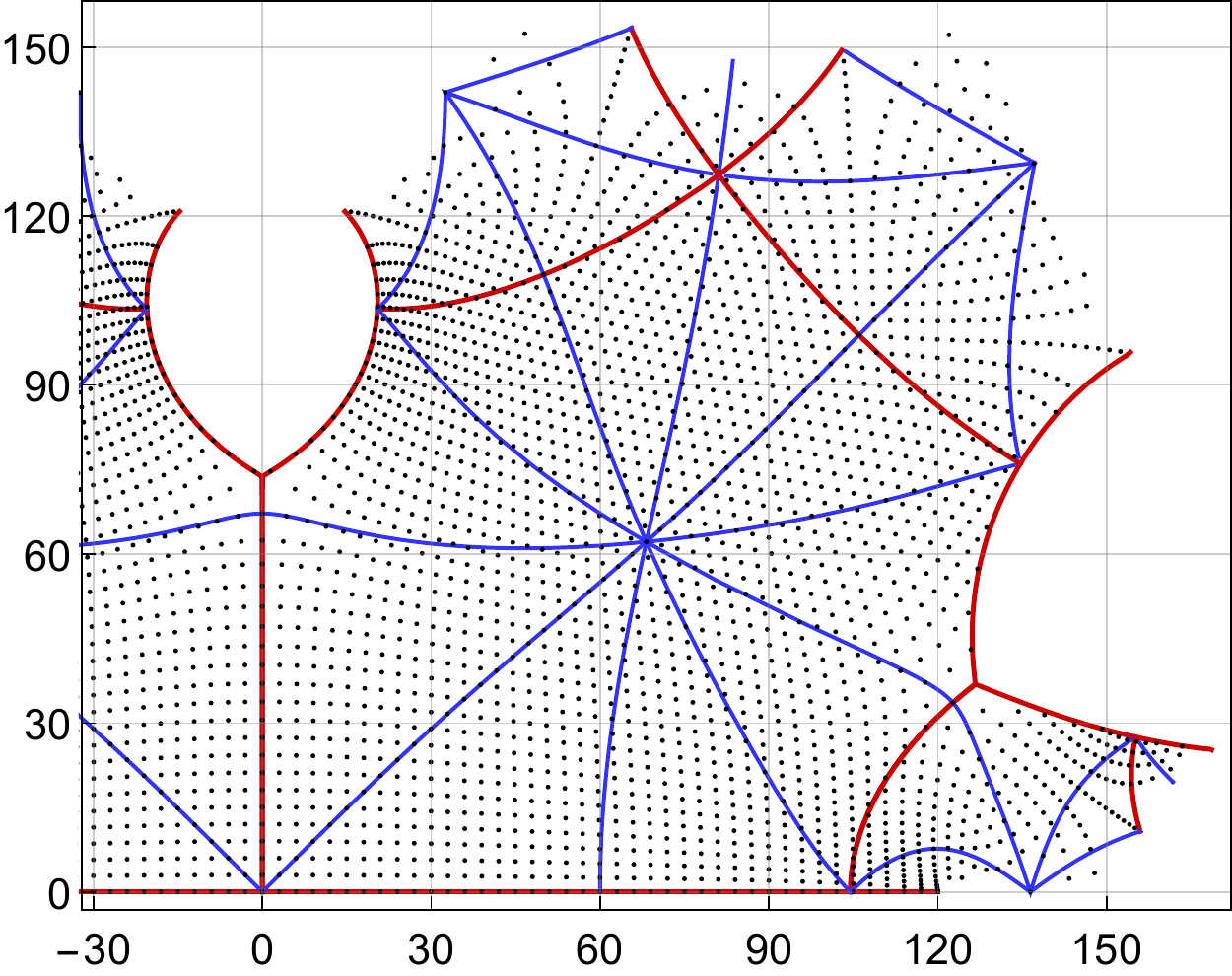}\hskip 13pt\hbox{}
  % in AngCoordsRot0.nb
    \end{center}
	\vspace{-5pt}
	\caption{The central image is the first quadrant of $R$ in $(\xi,\eta)$ coordinates, with the short-leg loops in red and the hypot-long-leg loops in blue (and the bounding blue loop $Q=0$ also dashed).  There is a dot for each pentagon in which both~$\xi$ and~$\eta$ are multiples of $3^\circ$.  The green ellipses are circles of radius $1/15$ under the metric $e_5$.  The inset above enlarges the quadrant\s top corner.  In the image below, the positions of the points have been adjusted, moving up from the $\xi$ axis, to give a conformal map of $(E_5,e_5)$, based on local isothermal coordinates $(u,v)$ [A7].}
	\label{fig:TractBoundaries}
\end{figure}

The ditri-pentagon loops cease to be geodesics after the rescaling; but no such issue arises for the tract-bounding loops, either short-leg or hypot-long-leg.  Those loops are geodesics both before and after the rescaling, since they are pointwise fixed by self-inverse, orientation-reversing isomorphisms.  If we work in $(x,y)$\_coordinates, each of those loops is an algebraic curve whose equation we can find via an appropriate Gr\"obner basis.  In Figure~\ref{fig:TractBoundaries}, the short-leg loops are drawn in red and the hypot-long-leg loops in blue.

Many of the regions in that figure are entire tracts --- triangles with red-blue-blue boundaries --- though their vertex angles are increasingly distorted out near the boundary of $R$.  When a ditri-pentagon loop collapses, however, each of the~$12$ tracts that it crosses collapses into a red-blue digon on one side of the collapse point and a blue-blue digon on the other.  The enlargement at the top of Figure~\ref{fig:TractBoundaries} shows one side of the collapse point at $(2\pi/3,\pm\pi)$.  There are three digons of each type on each sheet nestled below $(2\pi/3,+\pi)$: from left to right, a thin blue-blue digon, an entire tract, two thick red-blue digons, a tract, two thick blue-blue digons, a tract, and a thin red-blue digon.  Those $12$ digons are halves of the $12$ tracts that were crossed, the other halves being nestled above $(2\pi/3, -\pi)$.

\subsection{The hairy formulas}

To compute the first fundamental form of $(E_5,e_5)$ in $(\xi,\eta)$ coordinates, we set $x=\tan(\xi/2)$ and $y=\tan(\eta/2)$ and determine~$t_2$ and~$t_3$ from (\ref{eq:stformulas}).  Rotating the pentagon to make $\theta_3=0$, we set $\theta_4=2\arctan t_3$, $\theta_5=\theta_4+\eta$, $\theta_2=-2\arctan t_2$, and $\theta_1=\theta_2-\xi$.  Starting from (\ref{eq:BagOfEdges}), we then find [A5] that, for our unit-perimeter pentagons, 
\begin{equation*}
\label{eq:LineElement}
de_5^2=\mystrut{12pt}\mydepthstrut{6pt}\smash{\sum_{i<\mkern 1mu j}(d\theta\subj-d\theta_i)^2/25=E\, d\xi^2+2F\, d\xi\mkern 2mu d\eta+G\, d\eta^2}
\end{equation*}
where:
\begin{align}
\notag
{\textstyle E}&{\textstyle\null\!=\frac{75+225x^2+110y^2-15x^4+490x^2y^2-45y^4+27x^6-22x^4y^2+153x^2y^4-18x^6y^2+9x^4y^4+3x^6y^4}
		{50\,\mystrut{7.5pt}(15+22x^2+22y^2-9x^4+60x^2y^2-9y^4+6x^4y^2+6x^2y^4-x^4y^4)}}    \\[4 pt]
\label{eq:FFundForm}
{\textstyle F}&{\textstyle\null\!=\frac{xy\,(-55-50x^2-50y^2-27x^4+20x^2y^2-27y^4+6x^4y^2+6x^2y^4+x^4y^4)}
		{25\,\mystrut{7.5pt}(15+22x^2+22y^2-9x^4+60x^2y^2-9y^4+6x^4y^2+6x^2y^4-x^4y^4)}}  \\[4 pt]
\notag
{\textstyle
G}&{\textstyle\null\!=\frac{75+110x^2+225y^2-45x^4+490x^2y^2-15y^4+153x^4y^2-22x^2y^4+27y^6+9x^4y^4-18x^2y^6+3x^4y^6}
		{50\,\mystrut{7.5pt}(15+22x^2+22y^2-9x^4+60x^2y^2-9y^4+6x^4y^2+6x^2y^4-x^4y^4)}}
\end{align}

\begin{figure}
	\begin{center}
	\includegraphics[scale=0.7]{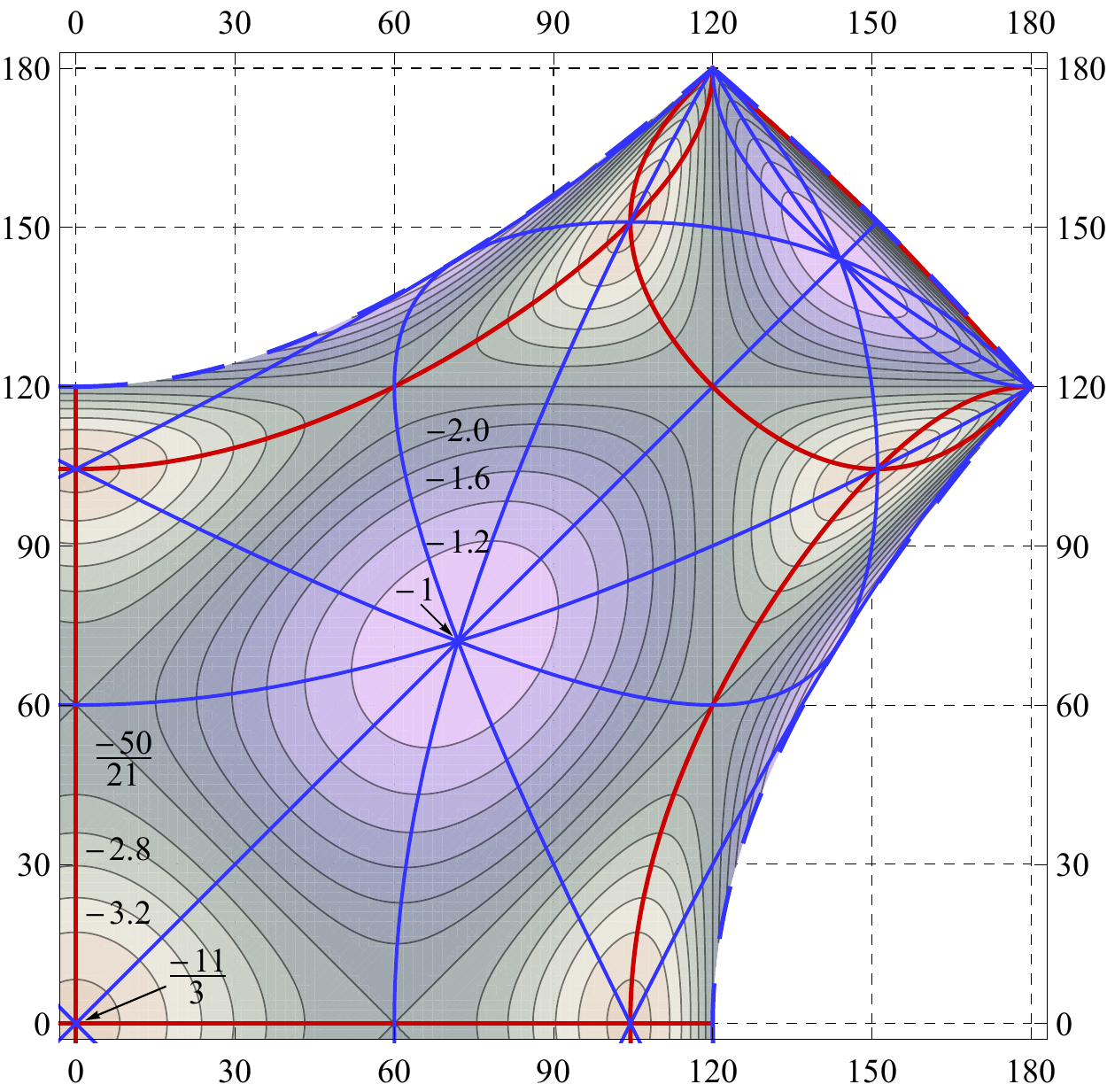}
	\end{center}
	\vspace{-10pt}
	\caption{The Gaussian curvature $K$ of $(E_5,e_5)$ over the first quadrant of $R$, as a function of $\xi$ and $\eta$.  The contour values run through the multiples of $0.2$ from $-1.2$ to $-3.6$, except that the one for $-2.4$ has been moved to $-50/21\approx -2.38$.  Intriguingly, $K$ is constant at that value along the ditri-pentagon loops, which  here plot as straight lines (including one from $(120^\circ,180^\circ)$ to $(180^\circ,120^\circ)$).}
	\label{fig:HomeDistrictK}
\end{figure}

We have written $E$, $F\mkern -1mu$, and $G$ in~(\ref{eq:FFundForm}) as rational functions of~$x$ and~$y$, with the radicand $Q$ from~(\ref{eq:subexpr}) in the denominators.  Written directly as functions of~$\xi$ and~$\eta$, they have many equivalent forms, none particularly simple.  The green ellipses in Figure~\ref{fig:TractBoundaries}, computed from~(\ref{eq:FFundForm}), show how nonisothermal the coordinates $(\xi,\eta)$ are.

The Brioschi formula\footnote{Francesco Brioschi (1824--1897) reworked the formula for calculating the curvature $K$ that Gauss had published in 1827 as the key step in proving his Theorema Egregium~\cite[p.~535]{Gray}.} then gives us the Gaussian curvature $K\mkern -1mu$, its formula~[A6] having $16$\_term factors of degree $12$.  The curvature varies from $-1$ at the $\frac{\pi}{5}$\_vertices down to $-11/3$ at the $\frac{\pi}{4}$\_vertices, as shown in Figure~\ref{fig:HomeDistrictK}.  The curvature is constant at $-50/21$ along the ditri-pentagon loops, and I don\t know why.  

\section{Initial conformal coordinates}
\label{sect:IsoCoords}  

A map of a Riemannian $2$\_manifold is conformal when it is the coordinate plane of an isothermal coordinate system, a system that solves the Beltrami equation.   As Spivak explains~\cite{SpivakIV}, Gauss showed how to solve the Beltrami locally for any $2$\_manifold whose metric is analytic (that is, real analytic).  Later work has shown that solutions exist locally even under much weaker assumptions; but Gauss\s technique suffices for us.  Given the complexity of~(\ref{eq:FFundForm}), though, a solution in closed form seems unlikely; so we approximate our isothermal coordinates numerically.

\subsection{Which isothermal coordinate system?}

Gauss constructs local isothermal coordinates $(u,v)$ and hence a local complex coordinate $w=u\pm v\imagunit$; but which one?  If $w_0=u_0+v_0\imagunit$ is one such, then any other $w_1=u_1+v_1\imagunit$ satisfies $w_1=h\circ w_0$ for some holomorphic function $h$ with nonzero derivative.  Gauss selects $h$ by requiring that all points along some line~$L$ through the origin of the $(\xi,\eta)$ plane retain, in the $(u,v)$ system, the same coordinates that they had in the $(\xi,\eta)$ system.  The coordinates of all other points are then adjusted, moving away from the line $L$, to make the $(u,v)$ system isothermal.  

The bottom image in Figure~\ref{fig:TractBoundaries} shows how this plays out if we choose $L$ to be the $\xi$\_axis.  We get a conformal map of $(E_5,e_5)$, with the proper tract angles of $\frac{\pi}{2}$, $\frac{\pi}{4}$, and~$\frac{\pi}{5}$.  That map, though, is drawn with local coordinates, and it peters out after just a few tracts.  As we move up from the $\xi$\_axis along either red side of the large decatract, we reach singularities of the type $z\mapsto z^{0.6667}$ --- presumably, actually $z\mapsto z^{2/3}$; and we are eventually unable to go any farther.

\begin{figure}
  \begin{center}
	\hskip 3 pt\includegraphics[scale=0.65]{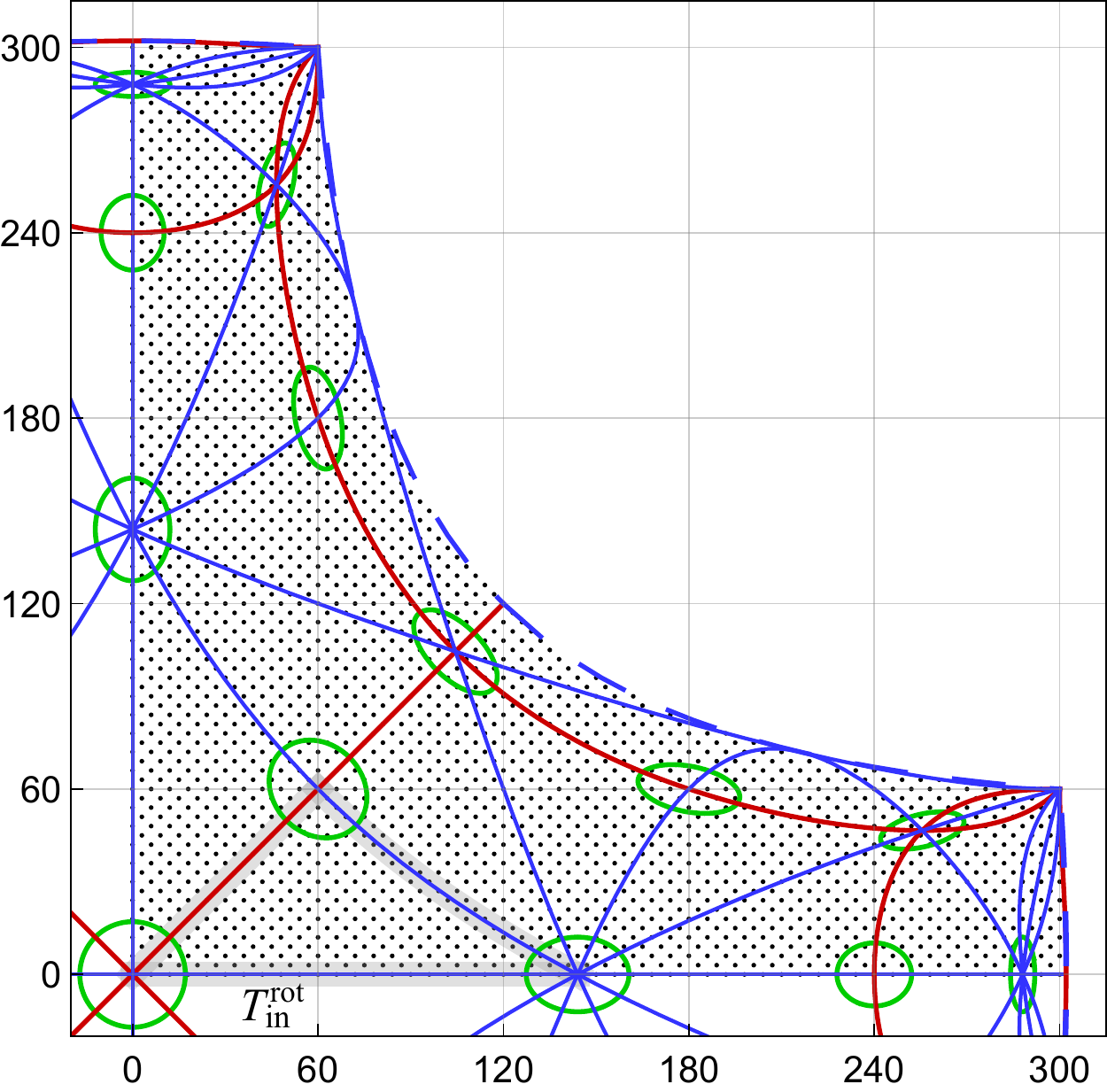} 
	         %from TractBoundaries.nb
  \end{center}
  \begin{center}
  \includegraphics[scale=0.642]{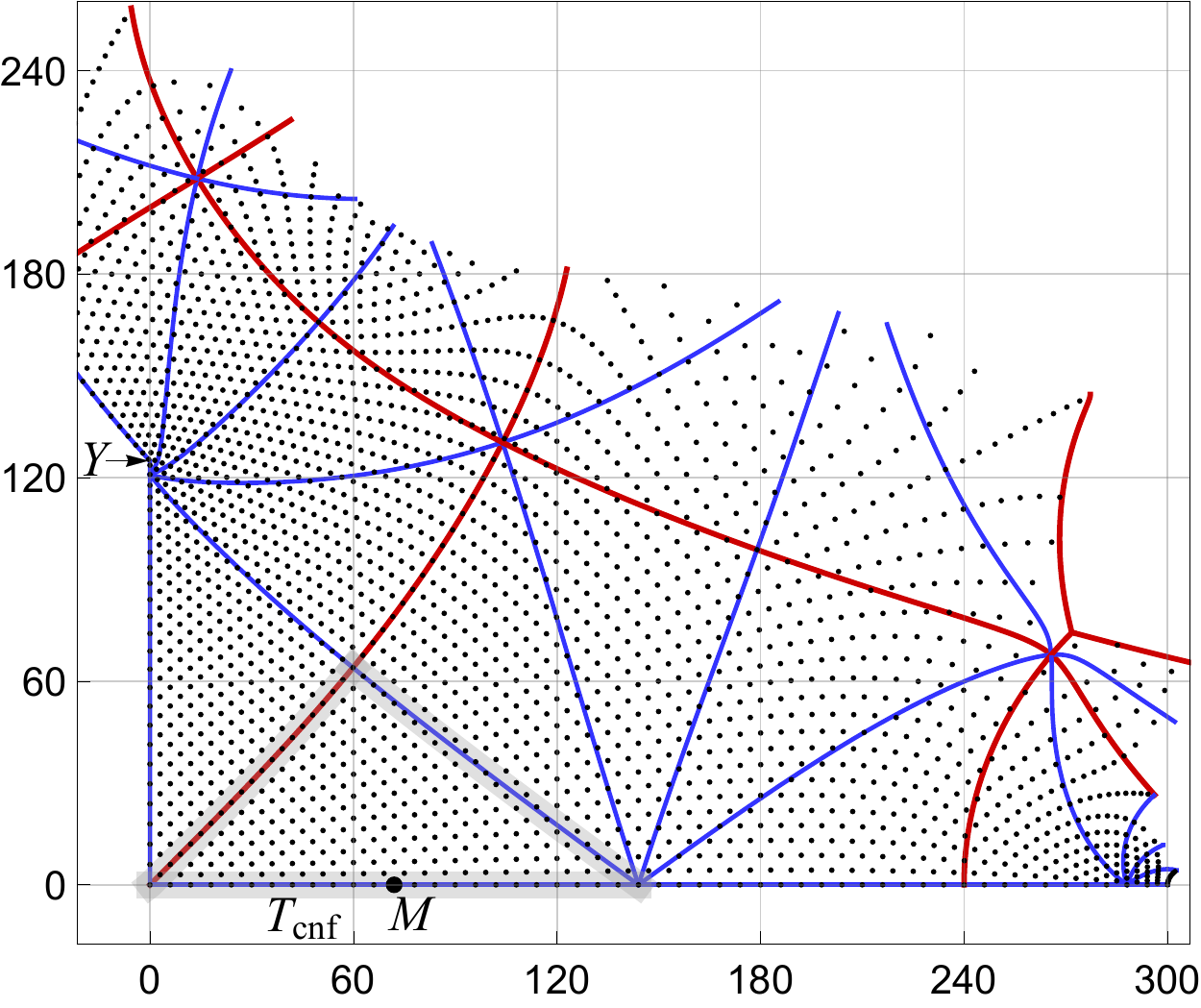}
  \end{center}
	\vspace{-5pt}
	\caption{The upper image shows the region $\eta\ge\lvert\xi\rvert$ of $R$, rotated clockwise $45^\circ$ to have $(\eta+\xi,\eta-\xi)$ as its coordinate axes.  The colors of the tract-bounding curves and the grid of dots are as in Figure~\ref{fig:TractBoundaries}.  In the lower image, the positions of the points have been adjusted, as we move northwest in $R$ from the main diagonal $\xi=\eta$, so as to make the new coordinates isothermal [A7].  The tract $T_\text{cnf}$ below is thus a conformal image of the tract that is depicted nonconformally as $T_\text{in}^\text{rot}$ above.  The singularity at $Y$ has exponent about $1.2348$.}
	\label{fig:Rot45}
\end{figure}

As another option, Figure~\ref{fig:Rot45} shows what happens if we pointwise fix the main diagonal, choosing as $L$ the line $\xi=\eta$.  The story is similar; but now, one of the first singularities that we reach has an exponent greater than $1$, so that angles there are expanded, rather than contracted.  That singularity is at $Y\mkern -2mu$, northwest of the origin along the line $\eta=-\xi$, and is of type roughly $z\mapsto z^{1.2348}$.

Our final conformal map, with its tracts adjusted to be precisely Poincar\'e projections of $(2,4,5)$\_triangles, won\t peter out.  Instead, it will fill out the entire Poincar\'e disk with infinitely many copies of the $240$ tracts of $E_5$.

\subsection{Gauss's technique}

For completeness, here is how Gauss computes isothermal coordinates, say when we pointwise fix the $\xi$-axis, as in Figure~\ref{fig:TractBoundaries}.  His key idea is to let the coordinate~$\xi$ take on complex values.

Gauss tells us to let $q(r)$ be some complex-valued function of a real variable $r$ that satisfies the ordinary differential equation
\begin{equation}
\label{eq:GaussODE}
q'(r)=\frac{-F-\imagunit\sqrt{EG-F^2}}{E}\bigl(q(r),r\bigr),
\end{equation}
evaluating $E$, $F$, and $G$ at $\xi=q(r)$ and $\eta=r$.  We set $w(\xi,\eta)$ to be $q(0)$ along the solution path of (\ref{eq:GaussODE}) that has $q(\eta)=\xi$, and hence passes through $(\xi,\eta)$.  This rule gives $w(\xi,0)=\xi$, and thus fixes the $\xi$\_axis pointwise, because $q(0)=\xi$ is the initial condition when $\eta=0$.  

Suppose that we move by an infinitesimal amount $(d\xi, d\eta)$ away from some point $(\xi,\eta)$, where the infinitesimals $d\xi$ and $d\eta$ satisfy
\[E\, d\xi+\bigl(F+\imagunit\sqrt{EG-F^2}\bigr)\,d\eta=0.\]
Since $q'(r)=d\xi/d\eta$, the vector $(d\xi, d\eta)$ is tangent to the solution curve of (\ref{eq:GaussODE}) that passes through $(\xi,\eta)$.  For an analytic metric, it follows~\cite[pp.~314--317\footnote{The $-\beta/\alpha$ in Equation (**) at the bottom of page 315 and top of page 316 should be $-\alpha/\beta$.}]{SpivakIV} that
\[ dw = c(\xi,\eta)\bigl(E\,d\xi+\bigl(F+\imagunit\sqrt{EG-F^2}\bigr)\,d\eta\bigl)
\]
for some smooth, complex-valued function $c(\xi,\eta)$.  Writing $c=a+b\imagunit$ and recalling that $w=u+v\imagunit$, we can calculate that
\begin{equation*}
\begin{pmatrix}u_\xi^2+v_\xi^2&u_\xi u_\eta+v_\xi v_\eta\\[2pt]
u_\xi u_\eta+v_\xi v_\eta&u_\eta^2+v_\eta^2\end{pmatrix}
=(a^2+b^2) E\begin{pmatrix}E&F\\F&G\end{pmatrix}\!.
\end{equation*}
Since $E$ is positive when $\xi$ and $\eta$ are real, the Euclidean metric of the $(u,v)$ coordinate plane thus pulls back to a metric on $E_5$ that is, as Gauss promised, a rescaling of the metric $e_5$.  

The square root in~(\ref{eq:GaussODE}) is an implementation challenge.   To approximate $w(\xi,\eta)$, we step along a solution curve of~(\ref{eq:GaussODE}), starting at $(\xi,\eta)$ and heading for the point~$\bigl(w(\xi,\eta), 0\bigr)$, visiting points $(\xi', \eta')$ in which $\eta'$ heads to zero through real numbers, but $\xi'$ is complex.  Along the way, the radicand $EG-F^2$ varies in~$\complexes$.  To avoid bad behavior, we must track its motion on the two-sheeted Riemann surface of the square-root function; a fixed branch cut won\t do.  For some starting points~$(\xi,\eta)$, however, the radicand passes through zero before the solution curve that we are following reaches $\eta'=0$, and bad behavior is then unavoidable.  In Figure~\ref{fig:TractBoundaries}, the two red boundaries of the large decatract, after we have moved up them a ways from the $\xi$-axis, consist entirely of such bad starting points.  That is what generates those $z\mapsto z^{2/3}$ branch points.

\section{Our conformal map}
\label{sect:ConformalMap}

For our final conformal map of $(E_5,e_5)$, it suffices, by symmetry, to conformally map any conformal image of a tract to any Poincar\'e projection of a $(2,4,5)$\_triangle.  If we choose our source conformal image and our target $(2,4,5)$\_triangle with care, the required conformal refinement will be almost linear, and we will be able to approximate it using an almost-linear polynomial.

\subsection{The setup}

The \Ltract\ $T_\text{cnf}$ in Figure~\ref{fig:Rot45} is a good starting point, and its straight hypotenuse is a plus.  For our source conformal image $H_\text{in}$, we shift $T_\text{cnf}$ leftward along the $(\eta+\xi)$\_axis to put the midpoint $M$ of its hypotenuse at the origin.  

\begin{figure}
  \begin{center}
	\raisebox{0.27 in}{\includegraphics[scale=0.2]{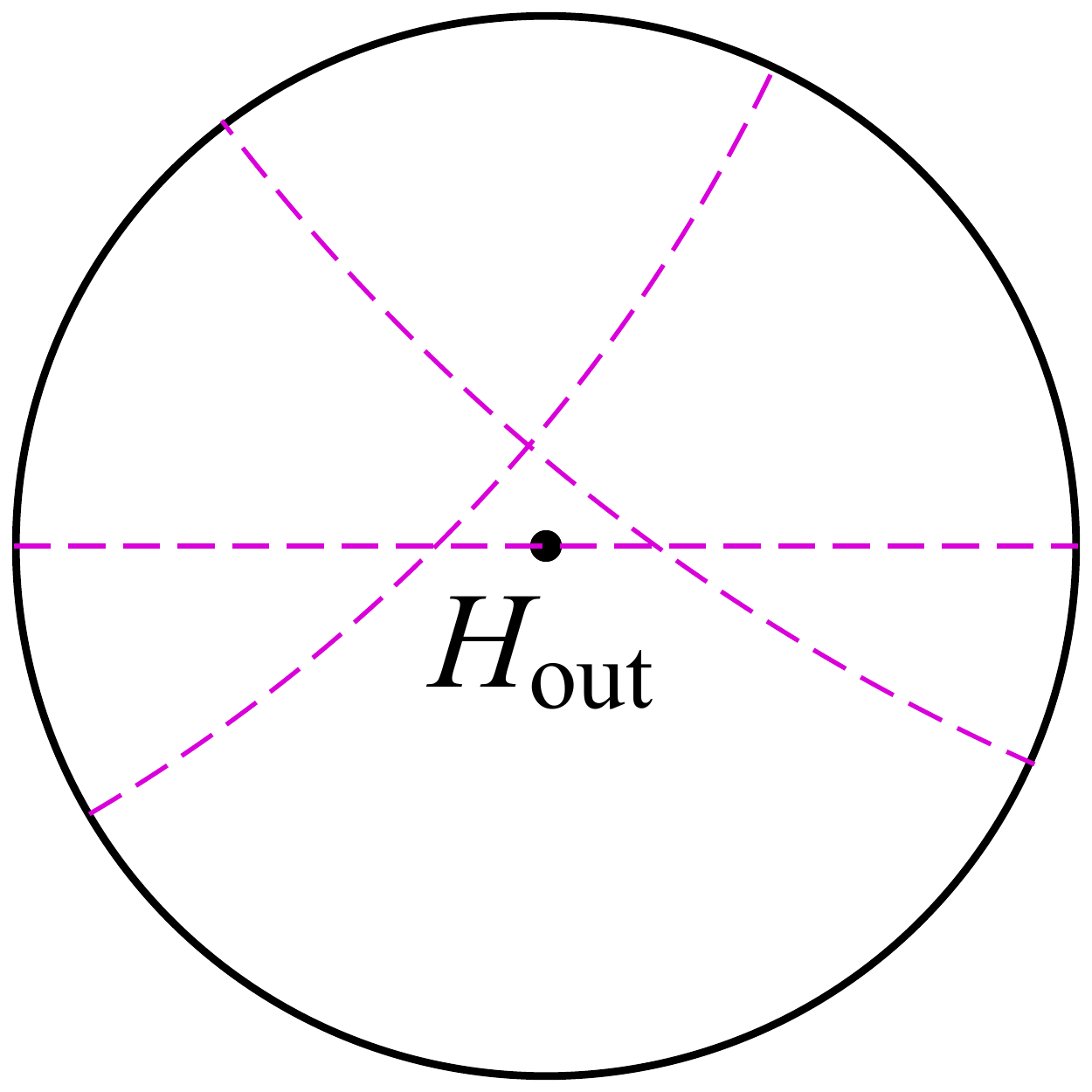}}\qquad
	\includegraphics[scale=0.6]{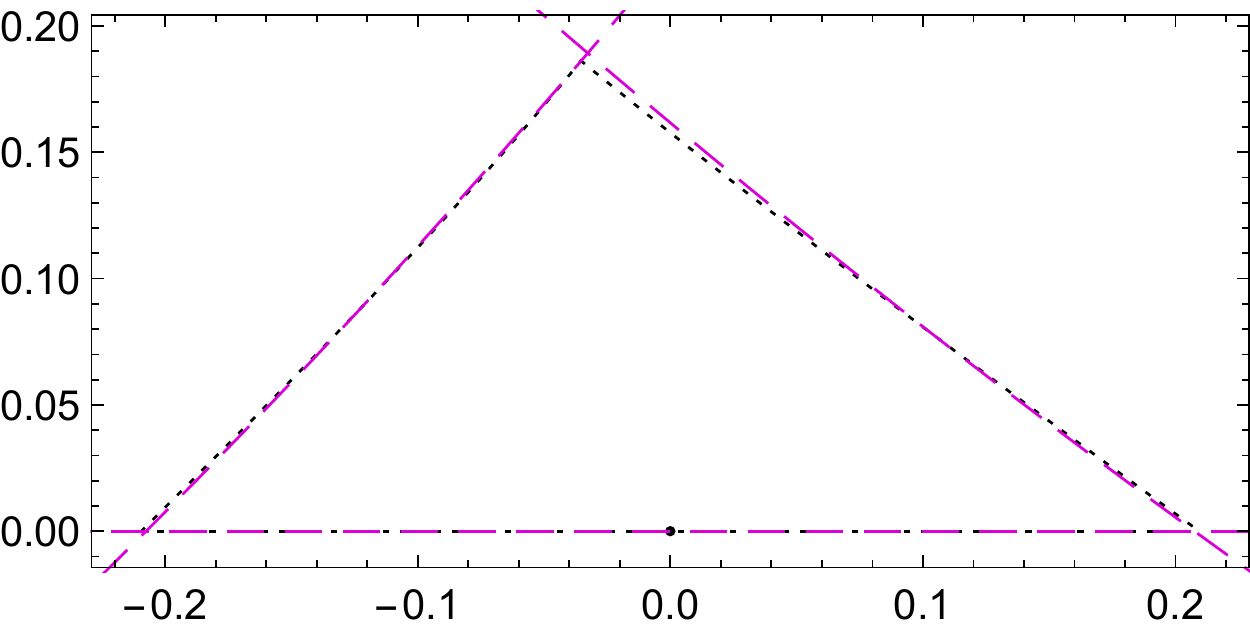}
  \end{center}
	\vspace{-10pt}
	\caption{Our target is the $(2,4,5)$\_triangle $H_\text{out}$ shown in dashed magenta, which has the midpoint of its hypotenuse at the origin.  The black dots on the right show the boundaries of $H_\text{in}/6$, their deviation peaking, near the upper right, at $0.004$.}
	\label{fig:TOut}
\end{figure}

As our target, we take the $(2,4,5)$\_triangle $H_\text{out}$ in Figure~\ref{fig:TOut}, which also has the midpoint of its hypotenuse at the origin.  There is a unique conformal bijection $h\colon H_\text{in}\to H_\text{out}$ that takes each vertex of $H_\text{in}$ to the vertex of $H_\text{out}$ with the same angle.  That map $h$ is holomorphic throughout $H_\text{in}$, including along its boundary and at its vertices, since the input and output vertex angles agree.

We approximate the map $h$ using a univariate polynomial $C$ of a complex variable whose coefficients are all real.  Any such polynomial preserves the real axis as a set, so we won\t bend the hypotenuse.  The right-hand side of Figure~\ref{fig:TOut} shows that the linear polynomial $C(z)=z/6$ already does pretty well.  (Using~$1/6$ as the linear coefficient here is a compromise;  it makes the dotted hypotenuse a skosh too long, but makes the dotted $\frac{\pi}{2}$\_vertex more than a skosh too low.)

We want~$C$ to be close to $h$ uniformly on $H_\text{in}$.  The max error will occur somewhere on the boundary; but we don\t know the point-to-point correspondence between the boundaries of $H_\text{in}$ and $H_\text{out}$ under~$h$, so there is no straightforward way to compute that max error.  We instead rate a candidate polynomial~$C$ using the max distance by which the short leg of the image $C(H_\text{in})$ deviates either inside or outside of the circular arc that forms the short leg of $H_\text{out}$, and similarly for their long legs.  The hypotenuse doesn\t deviate.  

Let $\dev(C)$ denote the maximum deviation, measured as a Euclidean distance orthogonal to the target circular arc.  Hyperbolic distance might be more natural, but would differ only by a factor close to $2$.  While $\dev(C)$ is only a lower bound on the true point-to-point max error, it seems likely to be in the right ballpark.  

\begin{figure}
  \begin{center}
    \includegraphics[scale=0.385]{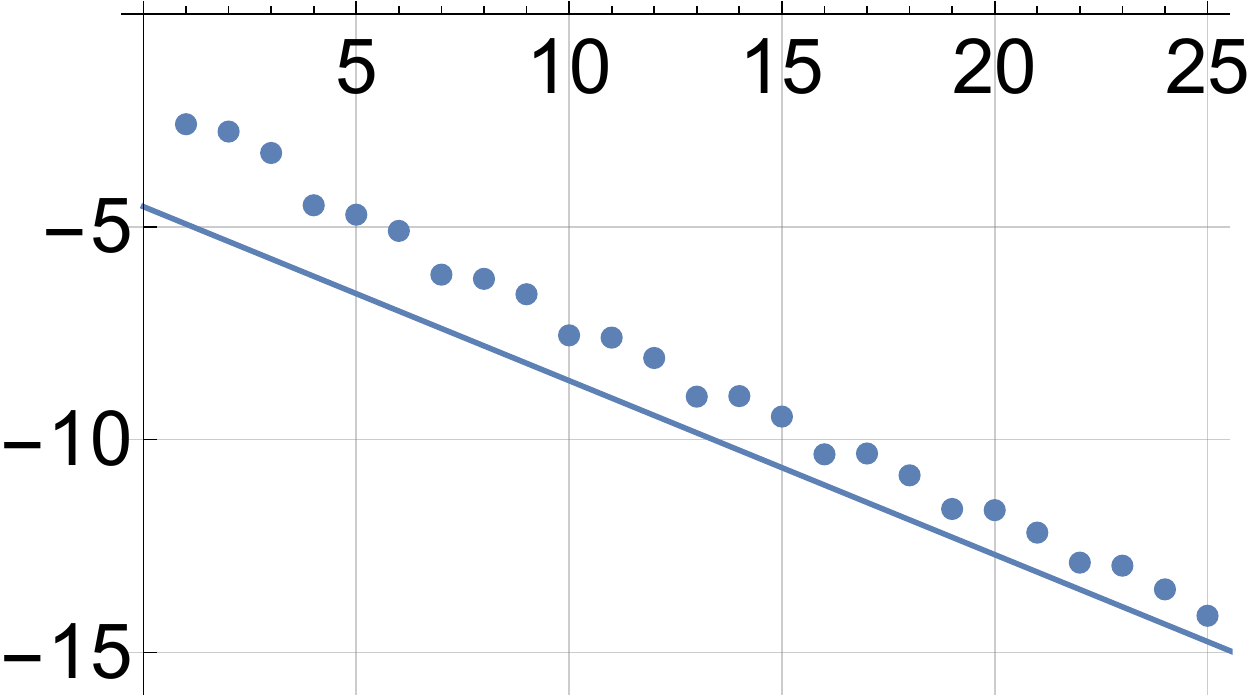}\quad 
	\includegraphics[scale=0.55]{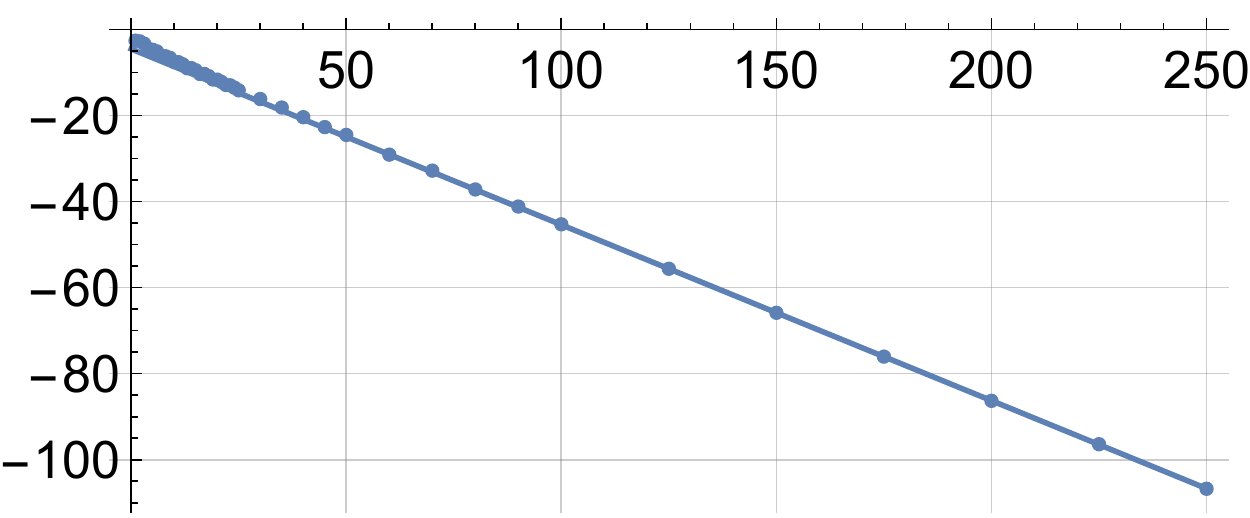}
  \end{center}
	\vspace{-10pt}
	\caption{A plot of $\log_{10}(\dev(C_d))$, with the apparent geometric trend line.}	
	\label{fig:ExpoConv}
\end{figure}

\subsection{The results}

For various degrees $d$ up to $250$, Mathematica looked for a polynomial $C_d$ of degree~$d$ that minimized the sum of the squared deviations at $2d+2$ sample points, $d+1$ of them Chebyshev distributed along each of the short and long legs~[A8].  Figure~\ref{fig:ExpoConv} is a semilog plot of $\dev(C_d)$.  The eventual convergence seems geometric, with each unit increase in $d$ reducing $\dev(C_d)$ by a factor of about $2.56$.  That isn\t surprising, since geometric convergence is common when approximating an analytic function with polynomials~\cite[Chap.~8]{Trefethen}.  The curves showing the deviation --- samples shown in [A8] --- come reasonably close to equioscillating.

To make our conformal map available to Mathematica users who want to plot pentagons on $E_5$, the appendix [A9] exports a function that computes the mapping from $T_\text{in}$ of Figure~\ref{fig:GlimpseAng} to $T_\text{out}$ of Figure~\ref{fig:Glimpse} with any requested accuracy up to $100$ digits, that limit being set by the accuracy of $C_{250}$.

\subsection{Canceling singularities}

As $d$ grows, our polynomials $C_d$ approach the Taylor series of the holomorphic function $h$ at the point $M$ in Figure~\ref{fig:Rot45}, the midpoint of the hypotenuse of~$T_\text{cnf}$.  The behavior of that Taylor series is determined largely by the singularities of $h$ that are close to $M\mkern -1mu$.  And those singularities are surprisingly easy to locate.

The lower image in Figure~\ref{fig:Rot45} is a complex-analytic chart for the Riemann surface $\RiemannEfive=(E_5^+,e_5^\smangle)$.  But that chart is only local, as the singularity $Y$ on the vertical axis in Figure~\ref{fig:Rot45} demonstrates.  The holomorphic function $h$ transforms that chart into another in which the tract $T_\text{cnf}$ has become precisely $H_\text{out}$, the Poincar\'e projection of a $(2,4,5)$\_triangle.  We could analytically continue that transformed chart by reflection across the sides of $H_\text{out}$ --- and then across the sides of those reflections --- to fill out the entire Poincar\'e disk with infinitely many copies of $\RiemannEfive$, each consisting of $240$ tracts.   So that transformed chart has no local-only restriction.  To enable this, the function $h$ must have singularities at $Y$ and at its conjugate~$\bar Y\mkern -1mu$, singularities that shrink angles by a factor of about $0.8099$, to undo their prior expansion by about $1.2348$.  And those will be the closest singularities of $h$ to~$M$, since the only justification for $h$ to have any singularity is to cancel out a singularity in the chart from Figure~\ref{fig:Rot45}.

Our polynomials~$C_d$ converge to~$h$ over~$T_\text{cnf}$ at a comfortable rate, captured in the slope of Figure~\ref{fig:ExpoConv}.  That happens, in part, because $Y$ is a goodly distance away from~$T_\text{cnf}$.  In addition, the phase of the vector $Y-M$, which is about $119.95^\circ$, is quite close to $2\pi/3$.  That helps to explain the oscillatory components of period $3$ that are visible, both in the max deviations graphed on the left in Figure~\ref{fig:ExpoConv} and in the coefficients of the individual polynomials $C_d$.  (The coefficients of $C_{50}$ are plotted in [A8].)

\subsection{The rescaling that normalizes curvature}

\begin{figure}
  \begin{center}
	\includegraphics[scale=0.5]{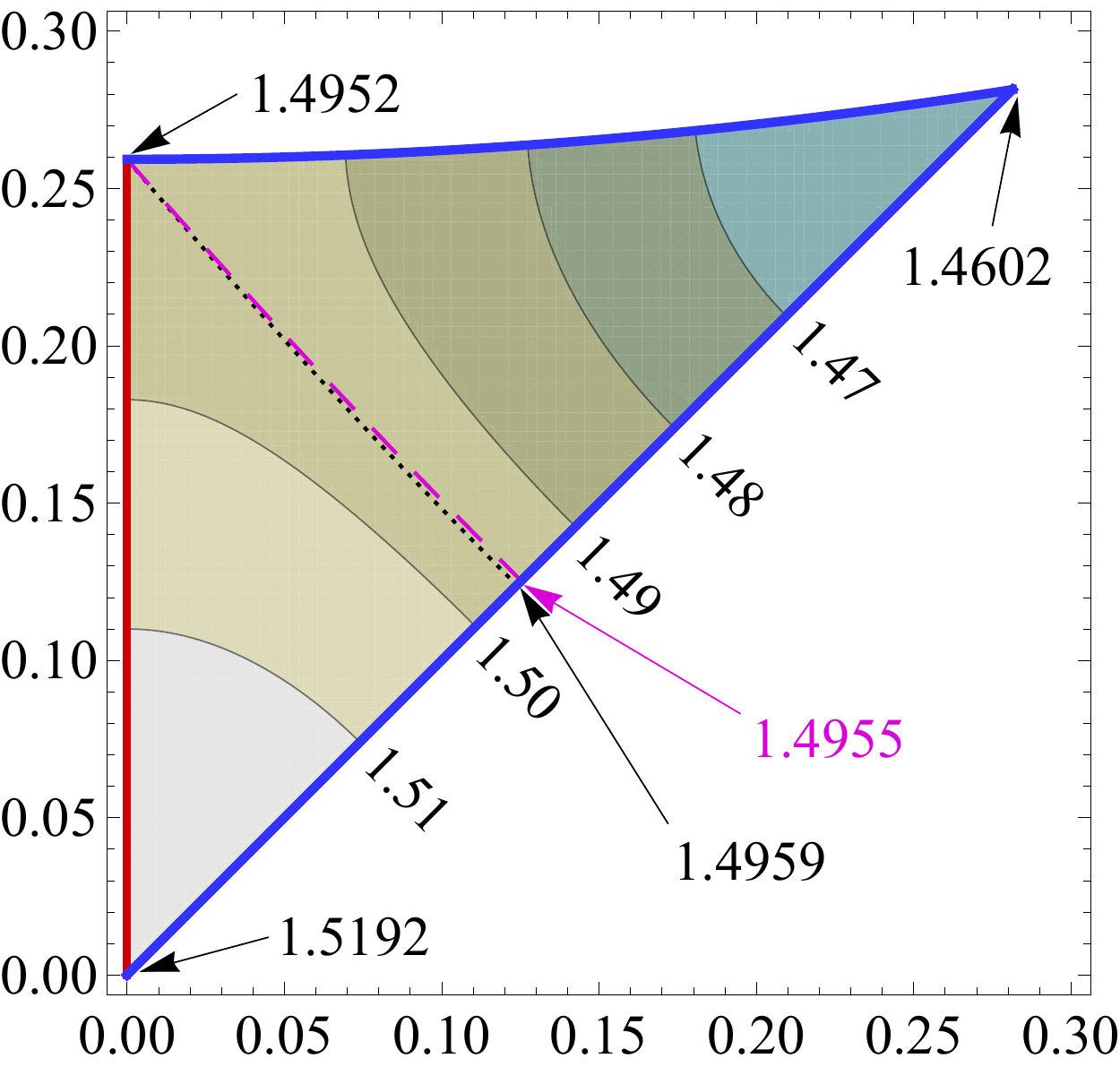}
  \end{center}
	\vspace{-10pt}
	\caption{The scale factor $s$ by which the manifold $(E_5,e_5)$ must have been stretched in order to make the curvature constant at $-1$, plotted over $T_\text{out}$.  The ditri-pentagon loop is black and dotted, while the close-by altitude loop is magenta and dashed.  The altitude loop stays a skosh farther away from the $\frac{\pi}{4}$\_vertex at the origin, so that it can travel through terrain that has been less stretched.}	
	\label{fig:ScaleFactor}
\end{figure}

Figure~\ref{fig:ScaleFactor} shows the rescaling $s$ that normalizes the curvature of $E_5$, computed by backing out, from our final conformal map, the effect of its Poincar\'e projection.

Before that rescaling, the curvature varies over $[-11/3\dd -1]$.  Multiplying all lengths by a constant $c$ divides all curvatures by $c^2$; so rescaling by any constant between $\sqrt{1}=1$ and $\sqrt{11/3}\approx 1.915$ would shift the interval of curvatures so that it contained $-1$.  As shown in Figure~\ref{fig:ScaleFactor}, the rescaling $s$ that makes the curvature constant at $-1$ uses a narrow band of scale factors near the middle of that range: the band $[1.4602\dd 1.5192]$, with a max-to-min ratio of only $1.04$.  

The rescaling $s$ has to behave the same on all tracts of $E_5$, like the Gaussian curvature $K$.  Thus, it isn\t surprising that $s$ varies roughly like $-K$ raised to a small power.  No formula of that type could be exact, however, since $K$ is constant along the ditri-pentagon loops, such as the one shown black and dotted in Figure~\ref{fig:ScaleFactor}, but $s$ is not.  The rescaling $s$ isn\t constant either along the altitude loops, such as the one shown magenta and dashed in Figure~\ref{fig:ScaleFactor}.

A child\s house gets plotted, on our conformal map, where the black-and-dotted ditri-pentagon loop crosses the diagonal.  That crossing point lies close to, but not on, the magenta-and-dashed altitude loop, which is the hyperbolic line joining the isosceles trapezoids at the two nearby $\frac{\pi}{2}$\_vertices, the distance between those loops as they cross the diagonal being about $0.48\%$ of the length of a tract hypotenuse.

\section{Acknowledgments}
Kudos to Gauss, who invented all of the differential geometry that we needed for this concrete cartographic problem, both concepts and methods.

The technical computing system Mathematica deserves a shout-out for its symbolic algebra, its accurate numerics, and its graphics primitives.  Its high-precision floating point uses \emph{significance arithmetic}, a cheaper cousin of interval arithmetic, to bound the round-off errors~\cite{Keiper}, and that was helpful in Section~\ref{sect:ConformalMap} [A8].

My gratitude goes also to Olivier B\'egassat, Fran\c cois Fillastre, John R.~Gilbert, Jean-Claude Hausmann, Michael Kapovich, Lance Ramshaw, James B.~Saxe, Jorge Stolfi, and Lloyd N.~Trefethen for helpful discussions.

\end{document}